# ASYMPTOTIC EQUIVALENCE OF NONPARAMETRIC AUTOREGRESSION AND NONPARAMETRIC REGRESSION

By Ion G. Grama and Michael H. Neumann

*Université de Bretagne-Sud and Friedrich-Schiller-Universität Jena*

It is proved that nonparametric autoregression is asymptotically equivalent in the sense of Le Cam's deficiency distance to nonparametric regression with random design as well as with regular nonrandom design.

**1. Introduction.** We assume that observations $X_0, \ldots, X_n$ from a stationary autoregressive process $(X_i)_{i=0,\ldots,n}$ are available which obey the model equation

$$(1) \qquad X_i = f(X_{i-1}) + \varepsilon_i, \qquad i = 1, \ldots, n,$$

where $(\varepsilon_i)_{i=1,\ldots,n}$ are i.i.d. random variables. The unknown autoregression function $f$ is then the target of statistical inference and the development of efficient estimators is a natural task for theoretically oriented statisticians. On the one hand, it has been recognized for a long time that commonly used estimators in model (1) have the same asymptotic behavior as corresponding estimators in nonparametric regression. A result of Robinson [26] concerns the pointwise equivalence of nonparametric kernel estimators and Neumann and Kreiss [22] extended this equivalence to the global behavior of nonparametric estimators. On the other hand, despite these well-known similarities between estimators, there is still a certain discrepancy in the current state of available theory in both contexts. While there is a very well developed asymptotic theory for optimal estimation in nonparametric regression, even up to the level of exact asymptotics (see, e.g., [13] or [24], for an overview), there is considerably less theory available in the case of nonparametric autoregression.

The purpose of the present paper is to bridge this gap between the two settings of nonparametric regression and autoregression by showing asymptotic equivalence on an abstract level. The theory of asymptotic equivalence of statistical experiments has been developed in Le Cam's [19] work. In the framework of nonparametric statistics, Brown and Low [4] proved that the Gaussian white noise experiment and nonparametric regression with nonrandom design and Gaussian errors are asymptotically equivalent in the sense that Le Cam's deficiency distance between them tends to zero. In [12, 14, 15, 23] the scope of asymptotic equivalence was extended to the nonparametric density estimation problem and to nonparametrically driven regression models. Moreover, asymptotic equivalence of nonparametric regression with random design and Gaussian white noise was shown in [2] while asymptotic equivalence of Poisson processes and Gaussian white noise was established in [3]. The issue of constructive asymptotic equivalence is considered in [25] and [5]. The asymptotic equivalence of a close relative of nonparametric autoregression, a diffusion experiment parametrized by the drift function, to Gaussian white noise experiments is proved in [8] and [7]. Milstein and Nussbaum [21] showed asymptotic equivalence of a nonparametric statistical model of small diffusion type and its discretization by a stochastic Euler difference scheme. These models deal with dependent observations in continuous time. However, asymptotic equivalence for models with dependent observations in discrete time where the noise is non-Gaussian seems to be a much more difficult issue.

In this paper we establish *local* equivalence of nonparametric autoregression (1) and nonparametric regression in the discrete-time setting. That is, the set of possible functions lies in a class $\Sigma_n(f_0)$ centered around some fixed function $f_0$ and shrinking in some appropriate norm as $n \to \infty$. Depending on additional prior smoothness assumptions on $f$, this class will nevertheless be rich enough for the transfer of minimax lower bounds from one to the other model. Under mild regularity assumptions stated below, the process $(X_i)_{i=0,\ldots,n}$ corresponding to $f_0$ has a stationary density $\psi_{f_0}$, say. We show asymptotic equivalence of the experiment given by (1) to nonparametric regression with random design as well as with regular nonrandom design. The former experiment corresponds to i.i.d. observations $(Y_1, \xi_1), \ldots, (Y_n, \xi_n)$ with

(2) $$Y_i = f(\xi_i) + \eta_i, \qquad i = 1, \ldots, n,$$

where $E(\eta_i|\xi_i) \equiv 0$. The basic assumption on the errors $\eta_i$ is that their Fisher information is the same as that of the $\varepsilon_i$'s. This includes the case of Gaussian errors as well as of errors having the same distribution as the $\varepsilon_i$. The $\xi_i$ are distributed according to the stationary density $\psi_{f_0}$ of the process corresponding to the central function $f_0$, regardless of the actual value of $f$.

We show also equivalence to nonparametric regression with regular nonrandom design, which corresponds to independent observations $Y_{n,1}, \ldots, Y_{n,n}$ obeying the model

$$(3) \qquad Y_{n,i} = f(t_{n,i}) + \eta_i, \qquad i = 1, \ldots, n,$$

where $E\eta_i = 0$. Here we will assume that the design points are regularly spaced with density $\psi_{f_0}$, that is, $\int_{-\infty}^{t_{n,i}} \psi_{f_0}(x) \, dx = (i - 1/2)/n$. We assume again that the Fisher information of $\eta_i$ is the same as the Fisher information of $\varepsilon_i$. Since Le Cam's equivalence relation is transitive, we also obtain as an immediate by-product asymptotic equivalence of nonparametric regression with random and regular nonrandom design. In the special case of Gaussian errors but under weaker smoothness assumptions on $f$, this equivalence also follows from the asymptotic equivalence of nonparametric regression with nonrandom design and Gaussian white noise [4] and the asymptotic equivalence of nonparametric regression with random design and Gaussian white noise [2].

At the end of Section 2 we discuss briefly how our results on asymptotic equivalence can be used to transfer well-known lower asymptotic bounds for the minimax risk in nonparametric regression to the case of nonparametric autoregression. Our local version of asymptotic equivalence does not allow an immediate transfer of upper asymptotic bounds; however, they could be independently proved by appeal to strong approximations of nonparametric estimators in both models (see [22] for details) or by direct computation of the risk of asymptotically optimal estimators.

**2. Assumptions and main results.** We start by introducing an appropriate functional parameter set. Consider the set of functions

$$\mathcal{F} = \left\{ f : \mathbb{R} \to \mathbb{R} : \sup_{x \in \mathbb{R}} |f(x)| \leq M \right\},$$

where $M < \infty$ is a constant. For any constants $\beta > 0$ and $L > 0$, let $\mathcal{H} = \mathcal{H}(\beta, L)$ be a Hölder ball, that is, the set of functions $f : \mathbb{R} \to \mathbb{R}$ satisfying

$$|f| \leq L, \qquad |f^{\lfloor \beta \rfloor}(x) - f^{\lfloor \beta \rfloor}(y)| \leq L|x - y|^{\beta - \lfloor \beta \rfloor}, \qquad x, y \in \mathbb{R}.$$

Here $\lfloor \beta \rfloor$ denotes the largest integer strictly less than $\beta$. The set of functional parameters is defined as

$$\Sigma = \mathcal{F} \cap \mathcal{H}(\beta, L).$$

Let $X_0$ be a random variable on the probability space $(\Omega, \mathcal{A}, P)$. Assume that we observe a sequence $X_1, \ldots, X_n$ which obeys

$$(4) \qquad X_i = f(X_{i-1}) + \varepsilon_i, \qquad i = 1, \ldots, n,$$

where $\varepsilon_1, \ldots, \varepsilon_n$ are i.i.d. with a given density $p$ that is continuous and positive on $\mathbb{R}$ and the function $f \in \Sigma$ is assumed to be unknown. It is easy to see that, for any $f \in \Sigma$, $P(X_{i+1} \in B | X_i = x) \geq \mu(B)$ holds for all $B \in \mathcal{B}$ and $x \in \mathbb{R}$, where $\mu$ is some measure not depending on $f$ with $\mu(\mathbb{R}) = \mu_0 > 0$. From Theorem 2.4.1 in [10] it follows that the uniform mixing coefficients (see Section 6.2) decay geometrically and, therefore, there exists a stationary density which we shall denote $\psi_f$.

Throughout the paper we shall assume that the observations (4) satisfy the following assumption:

(A1) The random variable $X_0$ has the stationary density $\psi_f(\cdot)$, which implies that the sequence $(X_i)_{i=0,\ldots,n}$ is in the stationary regime.

Note that the stationary density $\psi_f(\cdot)$ satisfies $\int_B \psi_f(x)\,dx \geq \mu(B)$, for all $B \in \mathcal{B}$.

Before we can state our main results on the approximation of the nonparametric autoregressive model by a nonparametric regression model, we have to introduce the basic concepts of *asymptotic equivalence*. Let $\mathcal{E}_l^n = (\Omega_l^n, \mathcal{A}_l^n, \{P_{l,f}^n, f \in \Sigma'\})$, $l = 1, 2$, be two sequences of statistical experiments indexed by $f$ in a subset $\Sigma' \subset \Sigma$. *The deficiency* of $\mathcal{E}_1^n$ with respect to $\mathcal{E}_2^n$ is defined as

$$\delta(\mathcal{E}_1^n, \mathcal{E}_2^n) = \sup_L \inf_{\delta^{(1)}} \sup_{\delta^{(2)}} \sup_{f \in \Sigma'} |E_{1,f}^n L(f, \delta^{(1)}) - E_{2,f}^n L(f, \delta^{(2)})|,$$

where the first supremum is taken over all decision problems with loss function $L$ with $0 \leq L \leq 1$, and the minimax value of the maximum difference in risks over $f \in \Sigma'$ is computed over all randomized statistical procedures $\delta^{(l)}$ for $\mathcal{E}_l^n$, $l = 1, 2$. According to Theorem 2 on page 15 in [20], the deficiency distance can alternatively be written as

$$\delta(\mathcal{E}_1^n, \mathcal{E}_2^n) = \inf_M \sup_{f \in \Sigma'} \tfrac{1}{2} \| M \cdot P_{1,f}^n - P_{2,f}^n \|_{\mathrm{Var}},$$

where $\| \cdot \|_{\mathrm{Var}}$ denotes the total variation distance and the infimum is taken over all Markov kernels $M$ on $\Omega_1^n \times \mathcal{A}_2^n$. *Le Cam's pseudodistance* between $\mathcal{E}_1^n$ and $\mathcal{E}_2^n$ is

$$\Delta(\mathcal{E}_1^n, \mathcal{E}_2^n) = \max\{\delta(\mathcal{E}_1^n, \mathcal{E}_2^n), \delta(\mathcal{E}_2^n, \mathcal{E}_1^n)\}.$$

Following [4], we say that the sequences $\mathcal{E}_1^n$, $n = 1, 2, \ldots$, and $\mathcal{E}_2^n$, $n = 1, 2, \ldots$, are *asymptotically equivalent* if

$$\Delta(\mathcal{E}_1^n, \mathcal{E}_2^n) \to 0 \qquad \text{as } n \to \infty.$$

To formulate our results we also need to impose the following regularity assumptions on the density $p(\cdot)$ of the innovations:

(A2)     (i) The density $p$ is positive on $\mathbb{R}$.

(ii) The log-likelihood function $l_p(x) = \log p(x)$ has three derivatives and satisfies, for some $\epsilon > 0$,

$$\int_\mathbb{R} \sup_{|u| \leq \epsilon} l_p''(x+u)^2 p(x)\, dx < \infty, \qquad \sup_{x \in \mathbb{R}} |l_p'''(x)| \leq c_1 < \infty.$$

(iii) The score $l_p'(x) = p'(x)/p(x)$ satisfies, for some $\epsilon > 0$ and any $\lambda < \infty$,

$$\int_\mathbb{R} \sup_{|u| \leq \epsilon} |l_p'(x+u)|^\lambda p(x)\, dx < \infty.$$

Assumption (A2) mainly requires the existence of three derivatives of $p(\cdot)$ and of the absolute moments of the corresponding scores. These types of assumptions can be related to the so-called Cramér conditions (see [20], page 102). Assumption (A2) is used here just for the sake of simplifying the proofs, but it is clear that they could be relaxed substantially. We refer to [15] for a relevant exposition of sufficient assumptions in the case of nonparametric models with independent observations.

In the sequel $q(\cdot)$ denotes a positive density which satisfies the following assumptions:

(A3)    (i) The log-likelihood function $l_q(x) = \log q(x)$ has three derivatives and satisfies, for some $\epsilon > 0$,

$$\int_\mathbb{R} \sup_{|u| \leq \epsilon} l_q''(x+u)^2 q(x)\, dx < \infty, \qquad \sup_{x \in \mathbb{R}} |l_q'''(x)| \leq c_1 < \infty.$$

(ii) The score $l_q'(x) = q'(x)/q(x)$ satisfies, for some $\epsilon > 0$ and any $\lambda < \infty$,

$$\int_\mathbb{R} \sup_{|u| \leq \epsilon} |l_q'(x+u)|^\lambda q(x)\, dx < \infty.$$

(iii) The Fisher information corresponding to the density $q(\cdot)$ is the same as that corresponding to $p(\cdot)$, that is,

$$I = \int_\mathbb{R} l_p'(x)^2 p(x)\, dx = \int_\mathbb{R} l_q'(x)^2 q(x)\, dx.$$

We state *local* versions of asymptotic equivalence, that is, we additionally assume that $f$ lies in a shrinking (as $n \to \infty$) neighborhood of some central function $f_0$. To get a meaningful result, we have to choose this neighborhood large enough such that it can be reached with a probability tending to 1 by an appropriate preliminary estimator. We fix any $\beta > 5/2$ and define

$$(5) \qquad \gamma_n = c\left(\frac{\log n}{n}\right)^{\beta/(2\beta+1)}, \qquad \gamma_n' = c'\left(\frac{\log n}{n}\right)^{(\beta-1)/(2\beta+1)}.$$

Here $\gamma_n$ and $\gamma'_n$ are the rates at which the function $f$ and its derivative $f'$ can be estimated in the model (4) and in the corresponding regression models. For any $f_0 \in \Sigma$, introduce the neighborhood

$$\Sigma^n_{f_0} = \{f \in \Sigma : f(x) = f_0(x), x \notin [A, B],$$
$$\|f - f_0\|_\infty \leq \gamma_n, \|f' - f'_0\|_\infty \leq \gamma'_n\},$$

where $A < B$ are two constants.

Our main results are the following two theorems which state the local asymptotic equivalence of our nonparametric autoregressive model to a nonparametric regression with random and nonrandom designs. We start with the case of random design.

THEOREM 2.1. *Let $\mathcal{E}^n_{f_0} = (\mathbb{R}^n, \mathcal{B}^n, \{P^n_f, f \in \Sigma^n_{f_0}\})$ be the local experiment based on observations $X_i$, $i = 0, \ldots, n$, obeying (A1) with $f \in \Sigma^n_{f_0}$. Suppose that the density $p(\cdot)$ satisfies assumption (A2). Let $\mathcal{G}^n_{f_0} = (\mathbb{R}^n, \mathcal{B}^n, \{Q^n_f, f \in \Sigma^n_{f_0}\})$ be the nonparametric regression model in which we observe*

$$(6) \qquad Y_i = f(\xi_i) + \eta_i, \qquad i = 1, \ldots, n,$$

*where $\eta_1, \ldots, \eta_n$ are i.i.d. with density $q(\cdot)$ obeying (A3), $\xi_1, \ldots, \xi_n$ are i.i.d. with the common density $\psi_{f_0}(\cdot)$, independent of $\eta_1, \ldots, \eta_n$, and $f \in \Sigma^n_{f_0}$ is unknown. Then, for all $\beta > 5/2$, the sequences of experiments $\mathcal{E}^n_{f_0}$, $n = 1, 2, \ldots$, and $\mathcal{G}^n_{f_0}$, $n = 1, 2, \ldots$, are asymptotically equivalent uniformly in $f_0 \in \Sigma$:*

$$\sup_{f_0 \in \Sigma} \Delta(\mathcal{E}^n_{f_0}, \mathcal{G}^n_{f_0}) \to 0 \qquad \text{as } n \to \infty.$$

Our second local result states asymptotic equivalence to the regression model with nonrandom design.

THEOREM 2.2. *Let $\mathcal{E}^n_{f_0} = (\mathbb{R}^n, \mathcal{B}^n, \{P^n_f, f \in \Sigma^n_{f_0}\})$ be the local experiment based on observations $X_i$, $i = 0, \ldots, n$, obeying assumption (A1) with $f \in \Sigma^n_{f_0}$. Assume that the density $p(\cdot)$ satisfies assumption (A2). Let $\mathcal{G}^n_{f_0} = (\mathbb{R}^n, \mathcal{B}^n, \{Q^n_f, f \in \Sigma^n_{f_0}\})$ be the nonparametric regression model in which we observe*

$$(7) \qquad Y_{n,i} = f(t_{n,i}) + \eta_i, \qquad i = 1, \ldots, n,$$

*where $\eta_1, \ldots, \eta_n$ are i.i.d. with density $q(\cdot)$ obeying assumption (A3). Furthermore, $t_{n,1}, \ldots, t_{n,n}$ are nonrandom design points chosen according to the density $\psi_{f_0}(\cdot)$, that is, $(i - 1/2)/n = \int_{-\infty}^{t_{n,i}} \psi_{f_0}(x)\, dx$, $i = 1, \ldots, n$, and $f \in \Sigma^n_{f_0}$ is unknown. Then, for all $\beta > 5/2$, the sequences of experiments $\mathcal{E}^n_{f_0}$, $n = 1, 2, \ldots$, and $\mathcal{G}^n_{f_0}$, $n = 1, 2, \ldots$, are asymptotically equivalent uniformly in $f_0 \in \Sigma$:*

$$\sup_{f_0 \in \Sigma} \Delta(\mathcal{E}^n_{f_0}, \mathcal{G}^n_{f_0}) \to 0 \qquad \text{as } n \to \infty.$$

REMARK 1. As a by-product of our main results, we obtain also asymptotic equivalence of nonparametric regression with random and regular nonrandom design. However, since we used a construction of the likelihood ratios based on a Skorokhod embedding rather than a KMT construction, the rate for the approximation error between the likelihood ratios of both models is presumably not the best possible one. We conjecture that the constraint $\beta > 5/2$ that was imposed for proving asymptotic equivalence of nonparametric autoregression and nonparametric regression can be further relaxed for the case of asymptotic equivalence of nonparametric regression with random and regular nonrandom design. It follows from the results in [4] and [2] that in the special case of Gaussian errors this equivalence holds even for $\beta > 1/2$.

REMARK 2. Our results on asymptotic equivalence in the Le Cam sense of nonparametric regression and autoregression can be used to transfer existing lower asymptotic efficiency bounds (when the loss is measured in the supremum norm) in nonparametric regression to the case of nonparametric autoregression. Indeed, it can be seen from the calculations in [9], Section 5, that a shrinking neighborhood of size $O((\log n/n)^{\beta/(2\beta+1)})$ around some central function $f_0$ is large enough for generating the desired risk bound. Hence, we can actually deduce these lower asymptotic efficiency bounds in the cases $\beta > 5/2$ which are covered by our results.

Owing to the local character of our results (asymptotic equivalence is proved for shrinking neighborhoods of $f_0$), we cannot directly use them for transferring upper asymptotic risk bounds. However, such bounds can be easily derived by straightforward calculations or by using asymptotic equivalence results between nonparametric estimators in both settings as given by strong approximations in [22].

The possibility of transferring asymptotic efficiency bounds on the basis of the asymptotic equivalence of experiments has been already known for a long time. This principle was applied by Korostelev and Nussbaum [18] for deducing asymptotic minimax bounds in nonparametric density estimation from known results in signal estimation in Gaussian white noise. On the basis of local equivalence results, Drees [11] transferred available lower asymptotic risk bounds from the Gaussian white noise model to the case of estimating an extreme value index.

**3. Proofs of the main theorems.** In this section we shall prove Theorem 2.1. Theorem 2.2 can be derived in the same way.

Our method of estimating the Le Cam distance $\Delta(\mathcal{E}_{f_0}^n, \mathcal{G}_{f_0}^n)$ runs as follows. Let $X_0, \ldots, X_n$ be the observations obeying assumption (A1) with $f \in \Sigma_{f_0}^n$ and let $(Y_1, \xi_1), \ldots, (Y_n, \xi_n)$ be the observations defined in Theorem 2.1. Denote by $L_{f,f_0}^{1,n}$ and $L_{f,f_0}^{2,n}$ the likelihood ratio processes of the experiments $\mathcal{E}_{f_0}^n$ and $\mathcal{G}_{f_0}^n$, respectively,

$$L_{f,f_0}^{1,n} = \frac{\psi_f(X_0)}{\psi_{f_0}(X_0)} \prod_{i=1}^n \frac{p(X_i - f(X_{i-1}))}{p(X_i - f_0(X_{i-1}))}$$

and

$$L_{f,f_0}^{2,n} = \prod_{i=1}^n \frac{q(Y_i - f(\xi_i))}{q(Y_i - f_0(\xi_i))}.$$

According to Proposition 2.2 in [23] (see also [20], page 16, for a similar assertion in the parametric context), the deficiency distance can be estimated as

$$(8) \qquad \Delta(\mathcal{E}_{f_0}^n, \mathcal{G}_{f_0}^n) \leq \sup_{f \in \Sigma_{f_0}^n} E_{\mathbf{P}} |\widetilde{L}_{f,f_0}^{1,n} - \widetilde{L}_{f,f_0}^{2,n}|,$$

where $\widetilde{L}_{f,f_0}^{1,n}$ and $\widetilde{L}_{f,f_0}^{2,n}$ are arbitrary versions of the likelihood ratios $L_{f,f_0}^{1,n}$ and $L_{f,f_0}^{2,n}$ constructed on a common probability space $(\Omega, \mathcal{F}, \mathbf{P})$ and distributed according to the central measure $P_{f_0}$. The versions $\widetilde{L}_{f,f_0}^{1,n}$ and $\widetilde{L}_{f,f_0}^{2,n}$ will be constructed in such a way that the right-hand side of (8) tends to zero as $n \to \infty$. Since this will hardly cause any confusion, we drop the tildes in the notation of $\widetilde{L}_{f,f_0}^{1,n}$ and $\widetilde{L}_{f,f_0}^{2,n}$. With this agreement inequality (8) can be written as

$$(9) \qquad \Delta(\mathcal{E}_{f_0}^n, \mathcal{G}_{f_0}^n) \leq \sup_{f \in \Sigma_{f_0}^n} E_{f_0} |L_{f,f_0}^{1,n} - L_{f,f_0}^{2,n}|.$$

The subscript $f_0$ at the expectation indicates that the measure $\mathbf{P}$ corresponds to the central measure $P_{f_0}$.

First, we give a bound for the $L_1$-distance on the right-hand side of (9) in terms of the Hellinger distance:

$$(10) \quad \tfrac{1}{2} E_{f_0} |L_{f,f_0}^{1,n} - L_{f,f_0}^{2,n}| \leq H(P_f^n, Q_f^n) \equiv \sqrt{E_{f_0}(\sqrt{L_{f,f_0}^{1,n}} - \sqrt{L_{f,f_0}^{2,n}})^2},$$

where $L_{f,f_0}^{1,n}$ and $L_{f,f_0}^{2,n}$ mean the corresponding versions of the likelihood ratios. Here $H(P,Q)$ denotes the Hellinger distance between two probability measures $P$ and $Q$. Following an idea originating from [23] in the context of density estimation and from [14] in the context of regression with independent observations, we shall use an analogue of the following property of the Hellinger distance for product measures (see Lemma 2.17 in [27]):

$$(11) \qquad H^2\!\left(\bigotimes_{l=1}^{K_n} P^{(l)}, \bigotimes_{l=1}^{K_n} Q^{(l)}\right) \leq \sum_{l=1}^{K_n} H^2(P^{(l)}, Q^{(l)}),$$

where $P^{(l)}$ and $Q^{(l)}$ are the measures corresponding to certain disjoint blocks of observations and $K_n$ is a sequence satisfying $K_n \to \infty$ and $K_n/n \to 0$. The size of these blocks will be chosen small enough so that one can get reasonable estimates for $H^2(P^{(l)}, Q^{(l)})$. It is clear that the estimate (11) is essentially based on the product structure of the measures $\bigotimes_{l=1}^{K_n} P^{(l)}$ and $\bigotimes_{l=1}^{K_n} Q^{(l)}$, and in general does not directly apply to the case of dependent observations.

In the particular context of the dependent data under consideration, we proceed as follows. Set $K_n = [n^{1/6}]$. Split the set of indices $\{1, \ldots, n\}$ into $K_n$ blocks,

$$\mathcal{I}_l = \left\{i : (l-1)\frac{n}{K_n} < i \le l\frac{n}{K_n}\right\}, \qquad l = 1, \ldots, K_n.$$

Denote by $m_l$ the number of elements in the block $\mathcal{I}_l$, that is, $m_l = \#\mathcal{I}_l = O(n^{5/6})$. Let $i_l$ be the first element in the set $\mathcal{I}_l$. Furthermore, let $\mathcal{F}_0$ be the trivial $\sigma$-field and, for $1 \le l \le K_n$,

$$\mathcal{F}_l = \sigma(X_0, \ldots, X_{i_l-1}; (Y_1, \xi_1), \ldots, (Y_{i_l-1}, \xi_{i_l-1})).$$

The likelihood ratio corresponding to the observations $X_0, \ldots, X_n$ can be written as the product

$$L_{f,f_0}^{1,n} = \prod_{l=0}^{K_n} L_{f,f_0}^{1,(l)}, \qquad L_{f,f_0}^{1,(l)} = \prod_{i \in \mathcal{I}_l} \frac{p(X_i - f(X_{i-1}))}{p(X_i - f_0(X_{i-1}))}, \qquad 1 \le l \le K_n,$$

where $L_{f,f_0}^{1,(0)} = \psi_f(X_0)/\psi_{f_0}(X_0)$ and $L_{f,f_0}^{1,(l)}$ is the conditional (given $\mathcal{F}_l$) likelihood ratio generated by $(X_i : i \in \mathcal{I}_l)$. Analogously, in the case of a regression experiment with random design, we have that

$$L_{f,f_0}^{2,n} = \prod_{l=0}^{K_n} L_{f,f_0}^{2,(l)}, \qquad L_{f,f_0}^{2,(l)} = \prod_{i \in \mathcal{I}_l} \frac{q(Y_i - f(\xi_i))}{q(Y_i - f_0(\xi_i))}, \qquad 1 \le l \le K_n,$$

where $L_{f,f_0}^{2,(0)} = 1$. A generalization of (11) to our setting with dependent random variables is given by the following result.

LEMMA 3.1.
$$H^2(P_{f,f_0}^n, Q_{f,f_0}^n) \le \sum_{l=0}^{K_n} \operatorname{ess\,sup} E_{f_0}((\sqrt{L_{f,f_0}^{1,(l)}} - \sqrt{L_{f,f_0}^{2,(l)}})^2 | \mathcal{F}_l).$$

PROOF. The proof of this assertion is adapted from that of Lemma 2.17 in [27]. We rewrite the Hellinger distance as

$$\tfrac{1}{2}H^2(L_{f,f_0}^{1,n}, L_{f,f_0}^{2,n}) = \tfrac{1}{2}E_{f_0}(\sqrt{L_{f,f_0}^{1,n}} - \sqrt{L_{f,f_0}^{2,n}})^2 = 1 - E_{f_0}\prod_{l=0}^{K_n}\sqrt{L_{f,f_0}^{1,(l)} L_{f,f_0}^{2,(l)}}.$$

For the last term one easily deduces

$$E_{f_0} \prod_{l=0}^{K_n} \sqrt{L_{f,f_0}^{1,(l)} L_{f,f_0}^{2,(l)}}$$

$$= E_{f_0}\left[\left(\prod_{l=0}^{K_n-1} \sqrt{L_{f,f_0}^{1,(l)} L_{f,f_0}^{2,(l)}}\right) E_{f_0}(\sqrt{L_{f,f_0}^{1,(K_n)} L_{f,f_0}^{2,(K_n)}}|\mathcal{F}_{K_n})\right]$$

$$\geq E_{f_0}\left[\prod_{l=0}^{K_n-1} \sqrt{L_{f,f_0}^{1,(l)} L_{f,f_0}^{2,(l)}}\right] \operatorname{ess\,inf} E_{f_0}(\sqrt{L_{f,f_0}^{1,(K_n)} L_{f,f_0}^{2,(K_n)}}|\mathcal{F}_{K_n}).$$

Continuing in the same way we obtain

$$E_{f_0} \prod_{l=0}^{K_n} \sqrt{L_{f,f_0}^{1,(l)} L_{f,f_0}^{2,(l)}}$$

$$\geq E_{f_0}[\sqrt{L_{f,f_0}^{1,(0)} L_{f,f_0}^{2,(0)}}] \prod_{l=1}^{K_n} \operatorname{ess\,inf} E_{f_0}(\sqrt{L_{f,f_0}^{1,(l)} L_{f,f_0}^{2,(l)}}|\mathcal{F}_l)$$

$$\geq \prod_{l=0}^{K_n}(1 - \operatorname{ess\,sup} \tfrac{1}{2} E_{f_0}((\sqrt{L_{f,f_0}^{1,(l)}} - \sqrt{L_{f,f_0}^{2,(l)}})^2|\mathcal{F}_l)).$$

Using the inequality $1 - \prod(1-a_i) \leq \sum a_i$, which is true for all $0 \leq a_i \leq 1$, we obtain the assertion of the lemma. $\square$

Hence, we have an analogue of (11) for the case of dependent random variables. Separability, which is equivalent to independence of the factors in $\prod_{l=0}^{K_n} \sqrt{L_{f,f_0}^{1,(l)} L_{f,f_0}^{2,(l)}}$, is just achieved by transition to the "worst case" which is appropriately expressed by $\operatorname{ess\,sup} E_{f_0}((\sqrt{L_{f,f_0}^{1,(l)}} - \sqrt{L_{f,f_0}^{2,(l)}})^2|\mathcal{F}_l)$.

Note that, since $p$ is positive on $\mathbb{R}$ and $\sup_{f \in \Sigma_{f_0}^n} \|f\|_\infty \leq M < \infty$, the condition $\rho < 1$ of Lemma 6.1 below is satisfied with some $\rho$ depending on $p$ and $M$. Then Lemma 6.1 and assumption (A2) imply, as $n \to \infty$,

(12)
$$\sup_{f_0 \in \Sigma} \sup_{f \in \Sigma_{f_0}^n} E_{f_0}((\sqrt{L_{f,f_0}^{1,(0)}} - \sqrt{L_{f,f_0}^{2,(0)}})^2|\mathcal{F}_0)$$

$$= \sup_{f_0 \in \Sigma} \sup_{f \in \Sigma_{f_0}^n} E_{f_0}(\sqrt{\psi_f(X_0)/\psi_{f_0}(X_0)} - 1)^2 \to 0.$$

Now Theorem 2.1 follows from Lemma 3.1, (12) and from the following assertion which provides us with bounds for the conditional Hellinger distance.

PROPOSITION 3.1. *Suppose that assumptions (A1)–(A3) are satisfied. Then there exists a construction of the sequences $X_0,\ldots,X_n$ and $(Y_1,\xi_1)$, $\ldots$, $(Y_n,\xi_n)$ on a common probability space such that*

$$\max_{1\leq l\leq K_n} \sup_{f_0\in\Sigma} \sup_{f\in\Sigma^n_{f_0}} \operatorname{ess\,sup} E_{f_0}((\sqrt{L^{1,(l)}_{f,f_0}} - \sqrt{L^{2,(l)}_{f,f_0}})^2|\mathcal{F}_l) = o(K_n^{-1}).$$

The proof of this proposition is postponed to Section 4.

In the case of comparing nonparametric autoregression and regression with regular nonrandom design, we proceed analogously. We use the same splitting of the set of indices $\{1,\ldots,n\}$ into blocks $\mathcal{I}_1,\ldots,\mathcal{I}_{K_n}$ as above. The pairs $(Y_{n,1},t_{n,1}),\ldots,(Y_{n,n},t_{n,n})$ are rearranged in such a way that

$$\left|\int_{-\infty}^{t_{n,n_l+i-1}} \psi_{f_0}(x)\,dx - \frac{i-1/2}{m_l}\right| = O(m_l^{-1}),$$

for all $i \in \{1,\ldots,m_l\}$, $l \in \{1,\ldots,K_n\}$. Then we write the likelihood ratio as

$$L^{3,n}_{f,f_0} = \prod_{l=0}^{K_n} L^{3,(l)}_{f,f_0}, \qquad L^{3,(l)}_{f,f_0} = \prod_{i\in\mathcal{I}_l} \frac{q(Y_{n,i}-f(t_{n,i}))}{q(Y_{n,i}-f_0(t_{n,i}))}, \qquad 1\leq l\leq K_n,$$

where $L^{3,(0)}_{f,f_0} = 1$. Let $\mathcal{F}'_0$ be the trivial $\sigma$-field and, for $l=1,\ldots,K_n$,

$$\mathcal{F}'_l = \sigma(X_0,\ldots,X_{i_l-1};Y_{n,1},\ldots,Y_{n,i_l-1}).$$

Using the same arguments as in the proof of Proposition 3.1 we obtain the following assertion.

PROPOSITION 3.2. *Suppose that assumptions (A1)–(A3) are satisfied. Then there exists a construction of the sequences $X_0,\ldots,X_n$ and $Y_{n,1}$, $\ldots$, $Y_{n,n}$ on a common probability space such that, as $n\to\infty$,*

$$\max_{1\leq l\leq K_n} \sup_{f_0\in\Sigma} \sup_{f\in\Sigma^n_{f_0}} \operatorname{ess\,sup} E_{f_0}((\sqrt{L^{1,(l)}_{f,f_0}} - \sqrt{L^{3,(l)}_{f,f_0}})^2|\mathcal{F}'_l) = o(K_n^{-1}).$$

Theorem 2.2 follows from Lemma 3.1, (12) and Proposition 3.2.

## 4. Proofs of Propositions 3.1 and 3.2.

PROOF OF PROPOSITION 3.1. Let $X_0,\ldots,X_n$ and $(Y_1,\xi_1),\ldots,(Y_n,\xi_n)$ be the observations generated according to (4) and (6). According to Theorem 5.1, there is a construction of the sequences $X_0,\ldots,X_n$ and $(Y_1,\xi_1),\ldots,(Y_n,\xi_n)$ on a common probability space which are coupled in such a way that the assertion of Theorem 5.1 holds true. Without loss of generality, we can assume that the sequences $X_0, \ldots, X_n$ and $(Y_1, \xi_1), \ldots, (Y_n, \xi_n)$ are already constructed on the probability space $(\Omega, \mathcal{F}^n, P_{f_0})$ endowed with the central measure $P_{f_0}$.

Recall that $m_l = \#\mathcal{I}_l = O(n^{5/6})$ is the number of indices in the set $\mathcal{I}_l$ and that $K_n = [n^{1/6}]$ is the number of blocks. Set, for brevity, $g(x) = f(x) - f_0(x)$. Since $f \in \Sigma_{f_0}^n$, we have $\|g\|_\infty \leq \gamma_n$ and $\|g'\|_\infty \leq \gamma_n'$. Since $\sup_x |l_p'''(x)| \leq c_1$ [by assumption (A2)(i)] and $\gamma_n^3 m_l = o(K_n^{-1/2})$, we obtain by a Taylor series expansion that

$$\log L_{f,f_0}^{1,(l)} = \sum_{i \in \mathcal{I}_l} g(X_{i-1}) l_p'(\varepsilon_i) + \tfrac{1}{2} \sum_{i \in \mathcal{I}_l} g(X_{i-1})^2 l_p''(\varepsilon_i) + o(K_n^{-1/2})$$
(13)
$$= T_1^{1,(l)} + T_2^{1,(l)} + o(K_n^{-1/2})$$

and, in the same way,

$$\log L_{f,f_0}^{2,(l)} = \sum_{i \in \mathcal{I}_l} g(\xi_i) l_q'(\eta_i) + \tfrac{1}{2} \sum_{i \in \mathcal{I}_l} g(\xi_i)^2 l_q''(\eta_i) + o(K_n^{-1/2})$$
(14)
$$= T_1^{2,(l)} + T_2^{2,(l)} + o(K_n^{-1/2}).$$

We introduce the set $A_l = A_{l,1} \cap A_{l,2}$, where

$$A_{l,1} = \{|T_1^{1,(l)} - T_1^{2,(l)}| \leq c_1 (\gamma_n)^{1/4} (\gamma_n')^{3/4} m_l^{1/4} \log m_l\},$$
$$A_{l,2} = \{|T_2^{1,(l)} - T_2^{2,(l)}| \leq v_n K_n^{-1/2}\},$$

and $v_n \to 0$ sufficiently slowly. An appropriate choice of the sequence $v_n$ is described in the course of the proof of Lemma 4.1 below. We bound the Hellinger distance between the partial likelihoods $L_{f,f_0}^{1,(l)}$ and $L_{f,f_0}^{2,(l)}$ as

(15) $$E_{f_0}((\sqrt{L_{f,f_0}^{1,(l)}} - \sqrt{L_{f,f_0}^{2,(l)}})^2 | \mathcal{F}_l) \leq R_1 + R_2,$$

say, where

$$R_1 = E_{f_0}((\sqrt{L_{f,f_0}^{1,(l)}} - \sqrt{L_{f,f_0}^{2,(l)}})^2 I_{A_l} | \mathcal{F}_l), \qquad R_2 = E_{f_0}(2(L_{f,f_0}^{1,(l)} + L_{f,f_0}^{2,(l)}) I_{\overline{A_l}} | \mathcal{F}_l).$$

First we bound $R_1$. On the set $A_l$, we get

$$|\log L_{f,f_0}^{1,(l)} - \log L_{f,f_0}^{2,(l)}| = O((\gamma_n)^{1/4}(\gamma_n')^{3/4} m_l^{1/4} \log m_l) + o(K_n^{-1/2}).$$

Since $\beta > 5/2$, we have $(\gamma_n)^{1/4}(\gamma_n')^{3/4} m_l^{1/4} \log m_l = o(K_n^{-1/2})$, which in turn implies that

$$|\sqrt{L_{f,f_0}^{1,(l)}/L_{f,f_0}^{2,(l)}} - 1| = |\exp(\tfrac{1}{2} \log L_{f,f_0}^{1,(l)} - \tfrac{1}{2} \log L_{f,f_0}^{2,(l)}) - 1| = o(K_n^{-1/2}).$$

Taking into account that $E_{f_0}(L^{2,(l)}_{f,f_0}|\mathcal{F}_l) = 1$, we get

$$
(16) \quad \begin{aligned} R_1 &= E_{f_0}(L^{2,(l)}_{f,f_0}(\sqrt{L^{1,(l)}_{f,f_0}/L^{2,(l)}_{f,f_0}} - 1)^2 I_{A_l}|\mathcal{F}_l) \\ &= o(K_n^{-1} E_{f_0}(L^{2,(l)}_{f,f_0}|\mathcal{F}_l)) = o(K_n^{-1}). \end{aligned}
$$

Now we shall bound $R_2$. Set $B_l = \{\log L^{1,(l)}_{f,f_0} \leq 1\}$ and $C_l = \{\log L^{2,(l)}_{f,f_0} \leq 1\}$. Then

$$
(17) \quad \begin{aligned} R_2 &\leq E_{f_0}(2(L^{1,(l)}_{f,f_0} I_{B_l} + L^{2,(l)}_{f,f_0} I_{C_l}) I_{\overline{A_l}}|\mathcal{F}_l) \\ &\quad + E_{f_0}(2(L^{1,(l)}_{f,f_0} I_{\overline{B_l}} + L^{2,(l)}_{f,f_0} I_{\overline{C_l}})|\mathcal{F}_l) \\ &\leq 4e P_{f_0}(\overline{A_l}|\mathcal{F}_l) + 2 E_{f_0}(L^{1,(l)}_{f,f_0} I_{\overline{B_l}}|\mathcal{F}_l) + 2 E_{f_0}(L^{2,(l)}_{f,f_0} I_{\overline{C_l}}|\mathcal{F}_l). \end{aligned}
$$

We will prove that

$$
(18) \quad P_{f_0}(\overline{A_l}|\mathcal{F}_l) = o(K_n^{-1}), \qquad P_{f_0}\text{-a.s.},
$$

and that

$$
(19) \quad \begin{aligned} E_{f_0}(L^{1,(l)}_{f,f_0} I_{\overline{B_l}}|\mathcal{F}_l) &= o(K_n^{-1}), \\ E_{f_0}(L^{2,(l)}_{f,f_0} I_{\overline{C_l}}|\mathcal{F}_l) &= o(K_n^{-1}), \qquad P_{f_0}\text{-a.s.} \end{aligned}
$$

Then, in conjunction with (15)–(17), we obtain the desired bound

$$
E_{f_0}((\sqrt{L^{1,(l)}_{f,f_0}} - \sqrt{L^{2,(l)}_{f,f_0}})^2|\mathcal{F}_l) = o(K_n^{-1}).
$$

Hence, it remains to prove (18) and (19).

First we prove (18). By Theorem 5.1 (with some $\lambda$ large enough) we have that

$$
(20) \quad P_{f_0}(\overline{A}_{l,1}|\mathcal{F}_l) = O(m_l^{-\lambda}) = o(K_n^{-1}).
$$

To complete the proof of (18) we shall prove the following bound.

LEMMA 4.1.

$$
(21) \quad P_{f_0}(\overline{A}_{l,2}|\mathcal{F}_l) = o(K_n^{-1}).
$$

PROOF. We shall use the fact that the Markov chain $X_0, \ldots, X_n$ is $\phi$-mixing. Decompose the set $\mathcal{I}_l$ as $\mathcal{I}_l = \mathcal{I}_l^{(1)} \cup \mathcal{I}_l^{(2)}$ so that $\mathcal{I}_l^{(1)}$ contains the first $c_0 \log m_l$ elements of the set $\mathcal{I}_l$ and $\mathcal{I}_l^{(2)}$ the remaining ones, where the positive constant $c_0$ will be chosen below. Let $l \in \{1, \ldots, K_n\}$ and let $i_{l,2}$ be the first element of $\mathcal{I}_l^{(2)}$. According to Lemma 6.2 in Section 6.2, we can construct a version $\widetilde{X}_{i_l,2-1}$ of the r.v. $X_{i_l,2-1}$ on the same probability space, such that $\widetilde{X}_{i_l,2-1}$ is independent of $\mathcal{F}_l$ and

(22) $$P_{f_0}(\widetilde{X}_{i_l,2-1} \neq X_{i_l,2-1}|\mathcal{F}_l) \leq \phi(X_{i_l,2-1}, \mathcal{F}_l) \leq c\rho^{c_0 \log m_l},$$

for some large enough constant $c_0$ and for some $\rho < 1$. Having constructed $\widetilde{X}_{i-1}$ for some $i \in \mathcal{I}_l^{(2)}$, we define recursively a version $\widetilde{X}_i$ of the r.v. $X_i$ on the same probability space, such that $\widetilde{X}_i$ is independent of $\mathcal{F}_l$ and of $\widetilde{X}_{i-c_0 \log m_l}, \ldots, \widetilde{X}_{i_l}, \varepsilon_{i-c_0 \log m_l}, \ldots, \varepsilon_{i_l}$ and

$$P_{f_0}(\widetilde{X}_i \neq X_i|\mathcal{F}_l) = \rho^{c_0 \log m_l}.$$

Choosing $c_0$ large enough, $\widetilde{X}_{i_l-1}, \ldots, \widetilde{X}_{i_{l+1}-1}$ satisfy

(23) $$P_{f_0}(\widetilde{X}_i \neq X_i, \forall i \in \{i_l - 1, \ldots, i_{l+1} - 1\}|\mathcal{F}_l) = m_l \rho^{c_0 \log m_l}$$
$$= o(K_n^{-1}).$$

Denote

$$\widetilde{T}_2^{1,(l)} = \tfrac{1}{2} \sum_{i \in \mathcal{I}_l} g(\widetilde{X}_{i-1})^2 l_p''(\varepsilon_i) = \tfrac{1}{2} \sum_{i \in \mathcal{I}_l} \zeta_i,$$

where $\zeta_i = g(\widetilde{X}_{i-1})^2 l_p''(\varepsilon_i)$. Since $\widetilde{T}_2^{1,(l)}$ is a sum of $c_0 \log m_l$-dependent r.v.'s, using Chebyshev's inequality, we obtain

$$P_{f_0}(|\widetilde{T}_2^{1,(l)} - E_{f_0}(\widetilde{T}_2^{1,(l)}|\mathcal{F}_l)| > v_n K_n^{-1/2}|\mathcal{F}_l)$$
$$\leq K_n/v_n E_{f_0}((\widetilde{T}_2^{1,(l)} - E_{f_0}(\widetilde{T}_2^{1,(l)}|\mathcal{F}_l))^2|\mathcal{F}_l)$$
$$= O\left(\frac{K_n}{v_n^2} \sum_{i \in \mathcal{I}_l} \sum_{j \in \mathcal{I}_l} E_{f_0}((\zeta_i - E_{f_0}(\zeta_i|\mathcal{F}_l))(\zeta_j - E_{f_0}(\zeta_j|\mathcal{F}_l))|\mathcal{F}_l)\right)$$
$$= O(v_n^{-2} K_n \gamma_n^4 m_l \log m_l) = o(v_n^{-2} n^{-1/3}).$$

Choosing $v_n$ such that $v_n \to 0$ and $o(v_n^{-2} n^{-1/3}) = o(K_n^{-1})$ we get

(24) $$P_{f_0}(|\widetilde{T}_2^{1,(l)} - E_{f_0}(\widetilde{T}_2^{1,(l)}|\mathcal{F}_l)| > v_n K_n^{-1/2}|\mathcal{F}_l) = o(K_n^{-1}).$$

By similar arguments for sums of independent random variables, we can show that

(25) $$P_{f_0}(|T_2^{2,(l)} - E_{f_0} T_2^{2,(l)}| > v_n K_n^{-1/2}|\mathcal{F}_l) = o(K_n^{-1}).$$

Taking into account that $E_{f_0} l_p''(\varepsilon_i) = E_{f_0} l_q''(\eta_i) = I$ [by Assumption (A3)(iii)] we obtain

$$E_{f_0}(\widetilde{T}_2^{1,(l)}|\mathcal{F}_l) - E_{f_0}(T_2^{2,(l)})$$
$$= O(\gamma_n^2 \log m_l) + \frac{I}{2} \sum_{i \in \mathcal{I}_l^{(2)}} E_{f_0}(g(\widetilde{X}_{i-1})^2|\mathcal{F}_l) - \frac{I}{2} \sum_{i \in \mathcal{I}_l^{(2)}} E_{f_0} g(\xi_i)^2.$$

Since $X_{i-1}$ and $\xi_i$ have the same density $\psi_f$ we get $E_{f_0}[g(X_{i-1})^2] = E_{f_0}[g(\xi_i)^2]$, and thus

$$(26) \quad E_{f_0}(\widetilde{T}_2^{1,(l)}|\mathcal{F}_l) - E_{f_0}(T_2^{2,(l)}) = O(\gamma_n^2 \log m_l) = o(K_n^{-1}).$$

By (23)–(26) we get

$$P_{f_0}(|T_2^{1,(l)} - T_2^{2,(l)}| \geq v_n K_n^{-1/2}|\mathcal{F}_l)$$
$$\leq P_{f_0}(|\widetilde{T}_2^{1,(l)} - T_2^{2,(l)}| \geq v_n K_n^{-1/2}|\mathcal{F}_l)$$
$$+ P_{f_0}(\widetilde{X}_i \neq X_i \ \forall i \in \{i_l - 1, \ldots, i_{l+1} - 1\}|\mathcal{F}_l)$$
$$= o(K_n^{-1}),$$

which proves (21). $\square$

Now we prove (19). We give a proof for the first bound; the second one can be proved in the same way. Changing the probability measure we obtain that

$$E_{f_0}(L_{f,f_0}^{1,(l)} I_{\overline{B}_l}|\mathcal{F}_l) = P_f(\overline{B}_l|\mathcal{F}_l) = P_f(\log L_{f,f_0}^{1,(l)} > 1|\mathcal{F}_l).$$

We shall prove that

$$(27) \quad P_f(\log L_{f,f_0}^{1,(l)} > 1|\mathcal{F}_l) = o(K_n^{-1}).$$

Indeed, proceeding as in the proof of (24) and using the fact that $\varepsilon_i = X_i - f_0(X_{i-1}) = X_i - f(X_{i-1}) + o(\gamma_n)$ and assumption (A2)(ii), one gets

$$(28) \quad P_f(|T_2^{1,(l)} - E_f(T_2^{1,(l)}|\mathcal{F}_l)| \geq cK_n^{-1/2}|\mathcal{F}_l) = o(K_n^{-1}).$$

Since $E_f(T_2^{1,(l)}|\mathcal{F}_l) = o(1)$, we get from (13) and (28),

$$(29) \quad P_f(\log L_{f,f_0}^{1,(l)} > 1|\mathcal{F}_l) \leq o(K_n^{-1}) + P_f(T_1^{1,(l)} > \tfrac{1}{2}|\mathcal{F}_l).$$

If we prove that

$$(30) \quad P_f(T_1^{1,(l)} > \tfrac{1}{2}|\mathcal{F}_l) = o(K_n^{-1}),$$

then we get, in conjunction with (29), that (27) holds.

To prove (30) we use the exponential Chebyshev's inequality for martingales. Since $\beta > 5/2$, by (5), we have $\gamma_n = o(n^{-5/12-3\delta})$, for some $\delta > 0$ small enough. Recall that $m_l = O(n^{5/6})$ and $\|g\|_\infty \leq \gamma_n$. Assume first that $n^{-\delta}|l'_p(\varepsilon_i)| \leq const$. Using Lemma 6.3,

$$P_f(T_1^{1,(l)} > \tfrac{1}{2}|\mathcal{F}_l)$$
$$\leq e^{-n^\delta} E_f\left(\exp\left(2\sum_{i\in\mathcal{I}_l} n^\delta g(X_{i-1})l'_p(\varepsilon_i)\right)\bigg|\mathcal{F}_l\right)$$

$$= e^{-n^\delta} E_f \left( \prod_{i \in \mathcal{I}_l} E_f(\exp(2n^\delta g(X_{i-1}) l'_p(\varepsilon_i)) | X_{i-1}) | \mathcal{F}_l \right)$$

$$\leq e^{-n^\delta} \prod_{i \in \mathcal{I}_l} \exp(cn^{2\delta} \gamma_n^2 E_f l'_p(\varepsilon_i)^2)$$

$$\leq e^{-n^\delta} \exp(cn^{2\delta} \gamma_n^2 m_l E_f l'_p(\varepsilon_i)^2),$$

where $c$ is a constant. The latter implies (30). If $n^{-\delta} |l'_p(\varepsilon_i)| \leq const$ is not satisfied, we use the same arguments with truncated scores $\hat{l}_i = \bar{l}_i - E_f(\bar{l}_i | X_{i-1})$, $\bar{l}_i = l'_p(\varepsilon_i) \mathbf{1}(|l'_p(\varepsilon_i)| \leq n^\delta)$ instead of the true scores $l'_p(\varepsilon_i)$. The term with the difference $l'_p(\varepsilon_i) - \hat{l}_i$ is bounded easily as before, using Chebyshev's inequality, the fact that $\varepsilon_i = X_i - f_0(X_{i-1}) = X_i - f(X_{i-1}) + o(\gamma_n)$ and assumption (A2)(iii):

$$P_f \left( \sum_{i \in \mathcal{I}_l} (g(X_{i-1}) l'_p(\varepsilon_i) - \hat{l}_i) > \tfrac{1}{2} | \mathcal{F}_l \right)$$
$$= O(\gamma_n^2 m_l E_f l'_p(\varepsilon_i)^2 \mathbf{1}(|l'_p(\varepsilon_i)| \geq n^\delta)) = O(K_n^{-1}).$$

Using the same types of arguments for sums of independent random variables we obtain

$$E_{f_0}(L^{2,(l)}_{f,f_0} I_{\overline{C}_l} | \mathcal{F}_l) = P_f(\overline{C}_l | \mathcal{F}_l) = P_f(\log L^{2,(l)}_{f,f_0} > 1 | \mathcal{F}_l) = o(K_n^{-1}),$$

which completes the proof of the first bound in (19). □

PROOF OF PROPOSITION 3.2. This proof is analogous to that of Proposition 3.1 and requires only a few minor modifications. Analogously to (13) and (14), we use the Taylor expansion

(31)
$$\log L^{3,(l)}_{f,f_0} = \sum_{i \in \mathcal{I}_l} g(t_{n,i}) l'_q(\eta_i) + \tfrac{1}{2} \sum_{i \in \mathcal{I}_l} g(t_{n,i})^2 l''_q(\eta_i) + o(K_n^{-1/2})$$
$$= T^{3,(l)}_1 + T^{3,(l)}_2 + o(K_n^{-1/2}).$$

Similarly to the calculations in the proof of Proposition 3.1, the closeness of $T^{1,(l)}_1$ and $T^{3,(l)}_1$ follows from Theorem 5.2, while that of $T^{1,(l)}_2$ and $T^{3,(l)}_2$ follows in complete analogy to the derivation of (21). □

**5. A functional strong approximation result.** In the proof of the main results we use the following strong approximation theorem. It can be viewed as an analogue of the functional strong approximation result established in [16] for sums of independent random variables.

Let $X_i$, $i = 1, \ldots, n$, and $(Y_i, \xi_i)$, $i = 1, \ldots, n$, be defined according to (4) and (6), respectively. Let $f_0 \in \Sigma$ and $f \in \Sigma_{f_0}^n$. We set

$$S_f^{1,(l)} = \sum_{i \in \mathcal{I}_l} (f - f_0)(X_{i-1}) l_p'(\varepsilon_i), \qquad S_f^{2,(l)} = \sum_{i \in \mathcal{I}_l} (f - f_0)(\xi_i) l_q'(\eta_i).$$

THEOREM 5.1. *Suppose that assumptions (A1)–(A3) are satisfied. Let $\lambda > 1$ be a constant. Then there are versions of the random variables $X_0, \ldots, X_n$ and $(Y_1, \xi_1), \ldots, (Y_n, \xi_n)$ on a common probability space such that, for $1 \leq l \leq K_n$,*

$$\sup_{f_0 \in \Sigma} \sup_{f \in \Sigma_{f_0}^n} \operatorname{ess\,sup} P_{f_0}(|S_f^{1,(l)} - S_f^{2,(l)}| > c(\lambda) r_n | \mathcal{F}_l) = O(m_l^{-\lambda}),$$

*where $r_n = (\gamma_n)^{1/4} (\gamma_n')^{3/4} m_l^{1/4} \log m_l + m_l^{-\lambda}$ and $c(\lambda)$ is a constant depending only on $\lambda$.*

The proof of this functional approximation result is based on a truncated Haar series expansion of $f - f_0$ and Lemma 5.1 below which provides a strong approximation result for partial sums with respect to a system of dyadic subintervals of $[A, B]$.

Define, for $j \geq 0$ and $k = 0, \ldots, 2^j$, $s_{j,k} = A + k 2^{-j}(B - A)$, and

$$I_{j,k} = (s_{j,k-1}, s_{j,k}], \qquad k = 1, \ldots, 2^j.$$

The Haar basis functions are defined via indicators as

$$h_0 = (B - A)^{-1/2} \mathbf{1}_{I_{0,1}},$$
$$h_{j,k} = (B - A)^{-1/2} 2^{j/2} (\mathbf{1}_{I_{j+1,2k-1}} - \mathbf{1}_{I_{j+1,2k}}) \qquad (j \geq 0; k = 1, \ldots, 2^j).$$

With a choice of the finest scale of the expansion, $j^* = j^*(n)$, described at the end of the proof of Theorem 5.1, we obtain a truncated Haar series expansion of $g = f - f_0$ as

$$g(x) = c_0(g) h_0(x) + \sum_{j=0}^{j^*} \sum_{k=1}^{2^j} c_{j,k}(g) h_{j,k}(x) + r_{j^*}(x),$$

where $c_0(g) = \int_A^B g(t) h_0(t) \, dt$, $c_{j,k}(g) = \int_A^B g(t) h_{j,k}(t) \, dt$, and $r_{j^*}(x)$ is the residual term. This yields that

$$|S_f^{1,(l)} - S_f^{2,(l)}|$$
$$\leq |c_0(g)| \left| \sum_{i \in \mathcal{I}_l} h_0(X_{i-1}) l_p'(\varepsilon_i) - \sum_{i \in \mathcal{I}_l} h_0(\xi_{i-1}) l_q'(\eta_i) \right|$$

$$+ \sum_{j=0}^{j^*} \sum_{k=1}^{2^j} |c_{j,k}(g)| \left| \sum_{i \in \mathcal{I}_l} h_{j,k}(X_{i-1}) l'_p(\varepsilon_i) - \sum_{i \in \mathcal{I}_l} h_{j,k}(\xi_{i-1}) l'_q(\eta_i) \right|$$

$$+ \left| \sum_{i \in \mathcal{I}_l} r_{j^*}(X_{i-1}) l'_p(\varepsilon_i) - r_{j^*}(\xi_{i-1}) l'_p(\eta_i) \right|$$

$$\leq (B-A)^{-1/2} |c_0(g)| |Z_{0,1}^{1,(l)} - Z_{0,1}^{2,(l)}|$$

$$+ \sum_{j=0}^{j^*} 2^{j/2} \sum_{k=1}^{2^j} |c_{j,k}(g)| (|Z_{j+1,2k-1}^{1,(l)} - Z_{j+1,2k-1}^{2,(l)}| + |Z_{j+1,2k}^{1,(l)} - Z_{j+1,2k}^{2,(l)}|)$$

$$+ \left| \sum_{i \in \mathcal{I}_l} r_{j^*}(X_{i-1}) l'_p(\varepsilon_i) - r_{j^*}(\xi_{i-1}) l'_p(\eta_i) \right|,$$

where

$$Z_{j,k}^{1,(l)} = \sum_{i \in \mathcal{I}_l} I(X_{i-1} \in I_{j,k}) l'_p(\varepsilon_i), \qquad Z_{j,k}^{2,(l)} = \sum_{i \in \mathcal{I}_l} I(\xi_{i-1} \in I_{j,k}) l'_q(\eta_i).$$

While the approximation-theoretic calculations are rather straightforward, the strong approximation result will require a lengthy proof based on Skorokhod embedding techniques. Let $\mathcal{I}_n = \{(j,k) : 0 \leq j \leq j^*, k = 1, \ldots, 2^j\}$.

LEMMA 5.1. *Suppose that assumptions (A1)–(A3) are satisfied. Then there exists a construction of the random variables $X_0, \ldots, X_n$ and $(Y_1, \xi_1)$, ..., $(Y_n, \xi_n)$ on a common probability space such that, for $1 \leq l \leq K_n$,*

$$\inf_{f_0 \in \Sigma} \operatorname{ess\,sup} P_{f_0}(|Z_{j,k}^{1,(l)} - Z_{j,k}^{2,(l)}| \leq C_\lambda (m_l 2^{-j})^{1/4} \log m_l, \forall (j,k) \in \mathcal{I}_n | \mathcal{F}_l)$$

$$= 1 - O(m_l^{-\lambda}).$$

To formulate the next theorem, we define

$$S_f^{3,(l)} = \sum_{i \in \mathcal{I}_l} (f - f_0)(t_{n,i}) l'(\eta_i).$$

THEOREM 5.2. *Suppose that assumptions (A1)–(A3) are satisfied. Let $\lambda > 1$ be a constant. Then there are versions of the random variables $X_0, \ldots, X_n$ and $Y_1, \ldots, Y_n$ on a common probability space such that, for $1 \leq l \leq K_n$,*

$$\sup_{f_0 \in \Sigma} \sup_{f \in \Sigma_{f_0}^n} \operatorname{ess\,sup} P_{f_0}(|S_f^{1,(l)} - S_f^{3,(l)}| > c(\lambda) r_n | \mathcal{F}'_l) = O(m_l^{-\lambda}),$$

*where $r_n = (\gamma_n)^{1/4} (\gamma'_n)^{3/4} m_l^{1/4} \log m_l + m_l^{-\lambda}$ and $c(\lambda)$ is a constant depending only on $\lambda$.*

The assertion of this theorem is a consequence of the following lemma. Set
$$Z_{j,k}^{3,(l)} = \sum_{i \in \mathcal{I}_l} I(t_{n,i} \in I_{j,k}) l'(\eta_i).$$

LEMMA 5.2. *Suppose that assumptions (A1)–(A3) are satisfied. Then there exists a construction of the random variables $X_0, \ldots, X_n$ and $Y_{n,1}, \ldots, Y_{n,n}$ on a common probability space such that, for $1 \leq l \leq K_n$,*

$$\inf_{f_0 \in \Sigma} \operatorname{ess\,sup} P_{f_0}(|Z_{j,k}^{1,(l)} - Z_{j,k}^{3,(l)}| \leq C_\lambda (m_l 2^{-j})^{1/4} \log m_l, \forall (j,k) \in \mathcal{I}_n | \mathcal{F}'_l)$$
$$= 1 - O(m_l^{-\lambda}).$$

The proofs of Lemmas 5.1 and 5.2 make use of a multiscale version of the Skorokhod embedding and are similar to the construction in [22]. We postpone these proofs to Section 5.2. Now we shall give proofs of Theorems 5.1 and 5.2.

5.1. *Proofs of Theorems 5.1 and 5.2.* As already indicated, the proofs of the theorems split into an approximation-theoretic and a stochastic part. The following lemma contains the approximation-theoretic facts needed for the proofs of Theorems 5.1 and 5.2.

LEMMA 5.3. *Let $c_0(g)$ and $c_{j,k}(g)$ be the Haar coefficients of a function $g$ defined above. Then:*

(i) $|c_0(g)| \leq (B-A)^{1/2} \|g\|_\infty$,
(ii) $|c_{j,k}(g)| \leq \min\{(B-A)^{1/2} 2^{-j/2} \|g\|_\infty (B-A)^{3/2} 2^{-3j/2-2} \|g'\|_\infty\}$,
(iii) $\|g - (c_0(g)h_0 + \sum_{j=0}^{j^*} \sum_{k=1}^{2^j} c_{j,k}(g) h_{j,k})\|_\infty \leq (B-A) 2^{-j^*-2} \|g'\|_\infty$.

PROOF. Assertion (i) follows from
$$|c_0(g)| \leq \|g\|_\infty \int |h_0(t)| \, dt \leq (B-A)^{1/2} \|g\|_\infty.$$

Analogously, we obtain that
$$|c_{j,k}(g)| \leq \|g\|_\infty \int |h_{j,k}(t)| \, dt \leq (B-A)^{1/2} 2^{-j/2} \|g\|_\infty.$$

Furthermore, it follows that
$$|c_{j,k}(g)| \leq (B-A)^{-1/2} 2^{j/2} \left| \int_A^B g(t)(\mathbf{1}_{I_{j+1,2k-1}}(t) - \mathbf{1}_{I_{j+1,2k}}(t)) \, dt \right|$$
$$\leq (B-A)^{-1/2} 2^{j/2} \int_{I_{j+1,2k-1}} |g(t) - g(t + (B-A) 2^{-j-1})| \, dt$$
$$\leq (B-A)^{3/2} 2^{-3j/2-2} \|g'\|_\infty,$$

which yields (ii).

Finally, we obtain from (ii) that

$$\left| g(x) - \left( c_0(g)h_0(x) + \sum_{j=0}^{j^*}\sum_{k=1}^{2^j} c_{j,k}(g)h_{j,k}(x) \right) \right|$$

$$\leq \sum_{j=j^*+1}^{\infty}\sum_{k=1}^{2^j} |c_{j,k}(g)h_{j,k}(x)|$$

$$\leq \sum_{j=j^*+1}^{\infty} (B-A)^{3/2}2^{-3j/2-2}\|g'\|_\infty (B-A)^{-1/2}2^{j/2}$$

$$= (B-A)2^{-j^*-2}\|g'\|_\infty. \qquad \square$$

Now we are in a position to prove Theorems 5.1 and 5.2.

PROOF OF THEOREM 5.1. Define

$$S_{f,j^*}^{1,(l)} = \sum_{i\in\mathcal{I}_l}\left[ c_0(g)h_0(X_{i-1}) + \sum_{j=0}^{j^*}\sum_{k=1}^{2^j} c_{j,k}(g)h_{j,k}(X_{i-1}) \right] l'_p(\varepsilon_i),$$

$$S_{f,j^*}^{2,(l)} = \sum_{i\in\mathcal{I}_l}\left[ c_0(g)h_0(\xi_i) + \sum_{j=0}^{j^*}\sum_{k=1}^{2^j} c_{j,k}(g)h_{j,k}(\xi_i) \right] l'_q(\eta_i)$$

and

$$R_{f,j^*}^{i,(l)} = S_f^{i,(l)} - S_{f,j^*}^{i,(l)}, \qquad i=1,2.$$

Define the event

$$D_l := \{|Z_{j,k}^{1,(l)} - Z_{j,k}^{2,(l)}| \leq C_\lambda (m_l 2^{-j})^{1/4}\log m_l, \forall\,(j,k)\in\mathcal{I}_n\},$$

where $C_\lambda$ is a constant. By Lemma 5.1, $P_{f_0}(\overline{D_l}|\mathcal{F}_l) = O(m_l^{-\lambda})$ with some choice of $C_\lambda$. By (i) and (ii) of Lemma 5.3, on the set $D_l$ it holds that

$$|S_{f,j^*}^{1,(l)} - S_{f,j^*}^{2,(l)}|$$

$$\leq (B-A)^{-1/2}\Bigg( |c_0(g)||Z_{0,1}^{1,(l)} - Z_{0,1}^{2,(l)}|$$

$$+ \sum_{j=0}^{j^*} 2^{j/2}\sum_{k=1}^{2^j} |c_{j,k}(g)|(|Z_{j+1,2k-1}^{1,(l)} - Z_{j+1,2k-1}^{2,(l)}|$$

$$+ |Z_{j+1,2k}^{1,(l)} - Z_{j+1,2k}^{2,(l)}|)\Bigg)$$

$$\leq C_\lambda \|g\|_\infty m_l^{1/4} \log m_l$$
$$+ C_\lambda \sum_{j=0}^{j^*} 2^{3j/2} \min\{2^{-j/2}\|g\|_\infty, (B-A)2^{-3j/2-2}\|g'\|_\infty\}$$
$$\times (m_l 2^{-j})^{1/4} \log m_l$$
$$= O((\|g\|_\infty + \|g\|_\infty^{1/4}\|g'\|_\infty^{3/4})m_l^{1/4} \log m_l).$$

The latter proves that, with some constant $c(\lambda)$ depending on $\lambda$,

$$(32) \qquad P_{f_0}(|S_{f,j^*}^{1,(l)} - S_{f,j^*}^{2,(l)}| > c(\lambda)r'_n|\mathcal{F}_l) = O(m_l^{-\lambda}),$$

where $r'_n = (\gamma_n)^{1/4}(\gamma'_n)^{3/4}m_l^{1/4} \log m_l$. By (iii) of Lemma 5.3 it holds that

$$(33) \qquad \begin{aligned} P_{f_0}(|R_{f,j^*}^{1,(l)}| \geq m_l^{-\lambda}) &\leq m_l^{2\lambda} E_{f_0}(R_{f,j^*}^{1,(l)})^2 \\ &\leq m_l^{2\lambda}(B-A)2^{-j^*-2}\|g'\|_\infty \sum_{i\in\mathcal{I}_l} E_{f_0}(l'_p(\varepsilon_i))^2. \end{aligned}$$

Choosing the finest level $j^*(n) = c^* \log m_l$, with some $c^*$ large enough, we obtain that

$$(34) \qquad P_{f_0}(|R_{f,j^*}^{1,(l)}| \geq m_l^{-\lambda}) = O(m_l^{-\lambda}).$$

Since the above bounds are uniform in $f_0 \in \Sigma$, from (32)–(34) and a similar bound for $R_{f,j^*}^{2,(l)}$ we conclude the assertion. $\square$

Theorem 5.2 can be proved in a similar way.

5.2. *Proofs of Lemmas* 5.1 *and* 5.2. We prove Lemma 5.1 only for $l=1$, since the proof for $l>1$ is completely analogous. The proof of Lemma 5.2 then requires only some obvious modifications and therefore will not be described here. To simplify notation we drop the index $l$ in the following, that is, we write $Z_{j,k}^1$, $Z_{j,k}^2$, $m$ instead of $Z_{j,k}^{1,(l)}$, $Z_{j,k}^{2,(l)}$, $m_l$, respectively.

PROOF OF LEMMA 5.1. Conditional on $X_0$ (which represents the information contained in $\mathcal{F}_0$), we construct a pairing of $X_1,\ldots,X_m$ and $(Y_1,\xi_1),\ldots,(Y_m,\xi_m)$ such that

$$\inf_{f_0\in\Sigma} \operatorname{ess\,sup} P_{f_0}(|Z_{j,k}^1 - Z_{j,k}^2| \leq C_\lambda(m2^{-j})^{1/4}\log m, (j,k)\in\mathcal{I}_n|X_0)$$
$$= 1 - O(m^{-\lambda})$$

is satisfied. In the following, all estimates are to be understood to hold uniformly in $f_0 \in \Sigma$.

The pairing of the random variables of both models is organized by a simultaneous Skorokhod embedding of $Z_{j,k}^1$ and $Z_{j,k}^2$ in a common set of Wiener processes $W_{j,k}$ assigned to the intervals $I_{j,k}$. We describe this embedding in detail for the autoregressive process. The embedding of $l_q(\eta_i)$ from the regression model is completely analogous and will be briefly mentioned only. Then we draw conclusions for the rate of approximation of $Z_{j,k}^1$ by $Z_{j,k}^2$, which will conclude the proof. An embedding scheme like this has already been developed in [22], in a different context. In view of some modifications and since we intend to provide a self-contained paper, we give a full proof of this lemma.

Let $W_{j,k}$, $(j,k) \in \mathcal{I}_n$, be independent Wiener processes. Apart from the coarsest resolution scale which corresponds to $j = 0$, we use each of these processes only on a finite time interval $[0, T_{j,k}]$, where the particular (nonrandom) values of the $T_{j,k}$ will be specified in part (iv) below. For the time being it is only important to know that $T_{0,k} = \infty$.

(i) *Embedding of $l_p(\varepsilon_1)$ and construction of $X_1$.*

First we define $l_p(\varepsilon_1)$ by a Skorokhod embedding in the Wiener processes mentioned above. Since $l_p(\varepsilon_1)$ does not necessarily define $X_1$ uniquely, we have to use perhaps an additional randomization to get $X_1$.

Let $k_1$ be that random number with $X_0 \in I_{j^*, k_1}$. Now we are going to represent $l_p(\varepsilon_1)$ by increments of the Wiener processes, preferably by those of $W_{j^*, k_1}$. However, since we want to use $W_{j^*, k_1}$ up to some prespecified time $T_{j^*, k_1}$ only, it might happen that this is not enough for representing $l_p(\varepsilon_1)$. In this case we additionally use a certain stretch of the process $W_{j^*-1, [k_1/2]}$, and so on. The Wiener processes which are potentially used for the representation of $l_p(\varepsilon_1)$ correspond to a containment relation of the dyadic intervals,

$$I_{j^*, k} \subseteq I_{j^*-1, [k/2]} \subseteq \cdots \subseteq I_{0, [k2^{-j^*}]},$$

where $[a]$ denotes the largest integer not greater than $a$. This means that we represent $l_p(\varepsilon_1)$ by the following Wiener process:

$$W^{(1)}(s) = \begin{cases} W_{j^*, k_1}(s), & \text{if } 0 \leq s \leq T_{j^*, k_1}, \\ W_{j^*, k_1}(T_{j^*, k_1}) + \cdots + W_{j+1, [k_1 2^{j+1-j^*}]}(T_{j+1, [k_1 2^{j+1-j^*}]}) \\ \quad + W_{j, [k_1 2^{j-j^*}]}\left(s - \sum_{l=j+1}^{j^*} T_{l, [k_1 2^{l-j^*}]}\right), \\ \qquad \text{if } T_{j^*, k_1} + \cdots + T_{j+1, [k_1 2^{j+1-j^*}]} \\ \qquad\qquad < s \leq T_{j^*, k_1} + \cdots + T_{j, [k_1 2^{j-j^*}]}. \end{cases}$$

($W^{(1)}$ is indeed a Wiener process on $[0, \infty)$, since $T_{0,k} = \infty$.)

According to Lemma A.2 of [17], there exists a stopping time $\tau^{(1)}$ such that the distribution of $W^{(1)}(\tau^{(1)})$ is equal to the conditional distribution of

$l_p(\varepsilon_1)$ given $X_0$. We define $\varepsilon_1$ in such a way that

$$l_p(\varepsilon_1) = W^{(1)}(\tau^{(1)}).$$

[This is achieved by first setting $l_p(\varepsilon_1)$ equal to $W^{(1)}(\tau^{(1)})$ and then defining $\varepsilon_1$ with the aid of an additional randomization according to its conditional distribution given $l_p(\varepsilon_1)$.] Finally, according to the model equation under $f_0$, we set $X_1 = f_0(X_0) + \varepsilon_1$.

To explain the following steps in a formally correct way, we introduce stopping times $\tau_{j,k}^{(i)}$, $i = 1, \ldots, m$, assigned to the corresponding Wiener process $W_{j,k}$. Define

$$\tau_{j,k}^{(0)} = 0, \qquad (j,k) \in \mathcal{I}_n.$$

To get $\tau_{j,k}^{(1)}$, we redefine all those $\tau_{j,k}^{(0)}$ which are assigned to Wiener processes $W_{j,k}$ that were used for representing $l_p(\varepsilon_1)$. According to the above construction we set

$$\tau_{j^*,k_1}^{(1)} = \tau^{(1)} \wedge T_{j^*,k_1}.$$

We redefine further

$$\tau_{j,[k_1 2^{j-j^*}]}^{(1)} = \begin{cases} [\tau^{(1)} - T_{j^*,k_1} - \cdots - T_{j+1,[k_1 2^{j+1-j^*}]}] \wedge T_{j,[k_1 2^{j-j^*}]}, \\ \quad \text{if } T_{j^*,k_1} + \cdots + T_{j+1,[k_1 2^{j+1-j^*}]} < \tau^{(1)}, \\ 0, \quad \text{otherwise.} \end{cases}$$

The remaining stopping times $\tau_{j,l}^{(1)}$ with $l \neq [k_1 2^{j-j^*}]$ keep their preceding values $\tau_{j,l}^{(0)} = 0$.

This procedure will be successively repeated for all other $\varepsilon_i$'s with the modification that we use only those parts of the Wiener processes which are still untouched by the previous construction steps.

(ii) *Embedding of $l_p(\varepsilon_i)$ and construction of $X_i$.*

Assume that $X_0, \ldots, X_{i-1}$ are already defined. Let $k_i$ be that random number with $X_{i-1} \in I_{j^*,k_i}$. Now we represent $l_p(\varepsilon_i)$ by parts of $W_{j^*,k_i}$, $W_{j^*-1,[k_i/2]}, \ldots, W_{0,[2^{-j^*}]}$, which have not been used so far.

First note that, because of the strong Markov property, these remaining increments $W_{j,[k_i 2^{j-j^*}]}(s + \tau_{j,[k_i 2^{j-j^*}]}^{(i-1)}) - W_{j,[k_i 2^{j-j^*}]}(\tau_{j,[k_i 2^{j-j^*}]}^{(i-1)})$ form independent Wiener processes, also independent of $X_0, \ldots, X_{i-1}$. Hence, gluing these parts together we obtain a Wiener process on $[0, \infty)$ which is independent of $X_0, \ldots, X_{i-1}$. This process is given as

$$W^{(i)}(s) = \begin{cases} W_{j^*,k_i}(s + \tau^{(i-1)}_{j^*,k_i}) - W_{j^*,k_i}(\tau^{(i-1)}_{j^*,k_i}), \\ \qquad \text{if } 0 \leq s \leq T_{j^*,k_i} - \tau^{(i-1)}_{j^*,k_i}, \\ \{W_{j^*,k_i}(T_{j^*,k_i}) - W_{j^*,k_i}(\tau^{(i-1)}_{j^*,k_i})\} + \cdots \\ \qquad + \{W_{j+1,[k_i 2^{j+1-j^*}]}(T_{j+1,[k_i 2^{j+1-j^*}]}) \\ \qquad - W_{j+1,[k_i 2^{j+1-j^*}]}(\tau^{(i-1)}_{j+1,[k_i 2^{j+1-j^*}]})\} \\ \qquad + \left\{ W_{j,[k_i 2^{j-j^*}]}\left(s - \sum_{l=j+1}^{j^*}(T_{l,[k_i 2^{l-j^*}]} - \tau^{(i-1)}_{l,[k_i 2^{j-l}]}) \right.\right. \\ \qquad \qquad \left. + \tau^{(i-1)}_{j,[k_i 2^{j-j^*}]}\right) \\ \qquad \qquad \left. - W_{j,[k_i 2^{j-j^*}]}(\tau^{(i-1)}_{j,[k_i 2^{j-j^*}]}) \right\}, \\ \qquad \text{if } \sum_{l=j+1}^{j^*}(T_{l,[k_i 2^{j-l}]} - \tau^{(i-1)}_{l,[k_i 2^{j-l}]}) < s \\ \qquad \qquad \leq \sum_{l=j}^{j^*}(T_{l,[k_i 2^{j-l}]} - \tau^{(i-1)}_{l,[k_i 2^{j-l}]}). \end{cases}$$

There exists a stopping time $\tau^{(i)}$ such that $W^{(i)}(\tau^{(i)})$ has the same distribution as the conditional distribution of $l_p(\varepsilon_i)$. We define $\varepsilon_i$ in such a way that

$$l_p(\varepsilon_i) = W^{(i)}(\tau^{(i)}),$$

and set $X_i = f_0(X_{i-1}) + \varepsilon_i$. [The definition of $\varepsilon_i$ is again achieved in two steps by first setting $l_p(\varepsilon_i)$ equal to $W^{(i)}(\tau^{(i)})$ and then defining $\varepsilon_i$ according to its conditional distribution.]

To complete this construction, it remains to define the stopping times $\tau^{(i)}_{j,k}$. These stopping times indicate up to which point the Wiener processes have been used in the first $i$ steps. Accordingly we set

$$\tau^{(i)}_{j,[k_i 2^{j-j^*}]} = \begin{cases} [\tau^{(i-1)}_{j,[k_i 2^{j-j^*}]} + (\tau^{(i)} - (T_{j^*,k_i} - \tau^{(i-1)}_{j^*,k_i}) - \cdots \\ \qquad - (T_{j+1,[k_i 2^{j+1-j^*}]} - \tau^{(i-1)}_{j+1,[k_i 2^{j+1-j^*}]}))] \\ \qquad \wedge T_{j,[k_i 2^{j-j^*}]}, \\ \qquad \text{if } (T_{j^*,k_i} - \tau^{(i-1)}_{j^*,k_i}) + \cdots \\ \qquad \qquad + (T_{j+1,[k_i 2^{j+1-j^*}]} \\ \qquad \qquad - \tau^{(i-1)}_{j+1,[k_i 2^{j+1-j^*}]}) < \tau^{(i)}, \\ \tau^{(i-1)}_{j,[k_i 2^{j-j^*}]}, \qquad \text{otherwise.} \end{cases}$$

For all $(j,l)$ with $l \neq [k_i 2^{j-j^*}]$ we define
$$\tau_{j,l}^{(i)} = \tau_{j,l}^{(i-1)}.$$

After embedding $l_p(\varepsilon_1), \ldots, l_p(\varepsilon_m)$ we arrive at stopping times $\tau_{j,k}^{(m)}$. The partial sums are connected to the Wiener processes by the relation

(35)
$$\begin{aligned} Z_{j,k}^1 &= \sum_{(u,v)\colon I_{u,v}\subseteq I_{j,k}} W_{u,v}(\tau_{u,v}^{(m)}) \\ &\quad + \sum_{i\colon 1\leq i\leq m, X_{i-1}\in I_{j,k}} \sum_{(u,v)\colon I_{u,v}\supset I_{j,k}} W_{u,v}(\tau_{u,v}^{(i)}) - W_{u,v}(\tau_{u,v}^{(i-1)}). \end{aligned}$$

(iii) *Embedding of $l_q(\eta_1), \ldots, l_q(\eta_m)$ and construction of $(Y_1, \xi_1), \ldots, (Y_m, \xi_m)$.*

This will be done in complete analogy to the construction described above. We define again stopping times $\widetilde{\tau}_{j,k}^{(i)}$ and obtain the following representation of the partial sums:

(36)
$$\begin{aligned} Z_{j,k}^2 &= \sum_{(u,v)\colon I_{u,v}\subseteq I_{j,k}} W_{u,v}(\widetilde{\tau}_{u,v}^{(m)}) \\ &\quad + \sum_{i\colon 1\leq i\leq m, Y_i\in I_{j,k}} \sum_{(u,v)\colon I_{u,v}\supset I_{j,k}} W_{u,v}(\widetilde{\tau}_{u,v}^{(i)}) - W_{u,v}(\widetilde{\tau}_{u,v}^{(i-1)}). \end{aligned}$$

(iv) *Choice of the values for $T_{j,k}$.*

To motivate our particular choice of the $T_{j,k}$ described below, we consider first two extreme cases. If $T_{j^*,k} = \infty$ for all $k$, then $Z_{j^*,k}^1$ and $Z_{j^*,k}^2$ are both completely represented by $W_{j^*,k}$. This leads indeed to a satisfactorily close connection of $Z_{j^*,k}^1$ and $Z_{j^*,k}^2$. On the other hand, this choice is unfavorable at scales $j \ll j^*$. Although we get immediately the upper estimate

$$|Z_{j,k}^1 - Z_{j,k}^2| \leq \sum_{l\colon I_{j^*,l}\subseteq I_{j,k}} |Z_{j^*,l}^1 - Z_{j^*,l}^2|,$$

the difference between $Z_{j,k}^1$ and $Z_{j,k}^2$ will be unnecessarily large. This is because, for $j \ll j^*$, $Z_{j,k}^1$ and $Z_{j,k}^2$ are then represented by too many different stretches of the Wiener processes $W_{j^*,l}$ with $I_{j^*,l} \subseteq I_{j,k}$.

On the other hand, if the $T_{j^*,k}$ are rather small, then $Z_{j^*,k}^1$ and $Z_{j^*,k}^2$ will be represented in large parts by stretches of Wiener processes $W_{u,v}$ which correspond to intervals $I_{u,v} \supset I_{j^*,k}$ with $j < j^*$. Once we are on a coarser scale $j < j^*$, we cannot guarantee that $Z_{j^*,k}^1$ and $Z_{j^*,k}^2$ are (mostly) generated by identical parts of the Wiener processes. Consequently, we would also get a suboptimal connection, this time for $Z_{j^*,k}^1$ and $Z_{j^*,k}^2$.

To find a good compromise between these two conflicting aims, we choose the $T_{j,k}$ as large as possible, however, with the additional property that, for $j \neq 0$, the stretches $[0, T_{j,k}]$ are used up with a high probability in the representation of both $l_p(\varepsilon_1), \ldots, l_p(\varepsilon_m)$ and $l_q(\eta_1), \ldots, l_q(\eta_m)$. Strictly speaking, we choose the $T_{j,k}$ in such a way that

$$
\text{(37)} \quad \operatorname*{ess\,sup} P_{f_0}\left(\sum_{i=1}^m \tau^{(i)} I(X_{i-1} \in I_{j,k}) < \sum_{(u,v)\,:\,I_{u,v} \subseteq I_{j,k}} T_{u,v}, \right.
$$
$$
\left. \forall\,(j,k) \in \mathcal{I}_n \setminus \{(0,k)\} \,\Big|\, X_0\right) = O(m^{-\lambda})
$$

and

$$
\text{(38)} \quad P_{f_0}\left(\sum_{i=1}^m \widetilde{\tau}^{(i)} I(Y_i \in I_{j,k}) < \sum_{(u,v)\,:\,I_{u,v} \subseteq I_{j,k}} T_{u,v},\, \forall\,(j,k) \in \mathcal{I}_n \setminus \{(0,k)\}\right)
$$
$$
= O(m^{-\lambda}).
$$

To this end, we study first the stochastic behavior of the above sums of stopping times assigned to the interval $I_{j,k}$.

Recall that the innovations $\varepsilon_i$ are assumed to be independent. According to the construction of the Skorokhod embedding described in [17], Appendix A.1, the randomness of $\tau^{(i)}$ is driven by some $U_i \sim Uniform[0, 1]$ from a sequence of independent random variables and by $\{W^{(i)}(s),\, 0 \leq s \leq \tau^{(i)}\}$. The vectors $(X_{i-1}, U_i)$ are of course also $\phi$-mixing as the $X_i$. Since, for $i \neq i'$, $\{W^{(i)}(s),\, 0 \leq s \leq \tau^{(i)}\}$ and $\{W^{(i')}(s),\, 0 \leq s \leq \tau^{(i')}\}$ are composed of disjoint stretches of the Wiener processes $W_{j,k}$ separated by stopping times, the random variables $\tau^{(i)} I(X_{i-1} \in I_{j,k})$ inherit the $\phi$-mixing property from the process $\{X_i\}$. Hence, we obtain by a Bernstein-type inequality for sums of $\phi$-mixing random variables (see, e.g., [10]) that

$$
\text{(39)} \quad \operatorname*{ess\,sup} P_{f_0}\left(\left|\sum_{i=1}^m \{\tau^{(i)} I(X_{i-1} \in I_{j,k}) - E[\tau^{(i)} I(X_{i-1} \in I_{j,k})]\}\right| \right.
$$
$$
\left. > C_\lambda \sqrt{m 2^{-j}} \log m \,\Big|\, X_0\right) = O(m^{-\lambda})
$$

and, analogously,

$$
\text{(40)} \quad P_{f_0}\left(\left|\sum_{i=1}^m \{\widetilde{\tau}^{(i)} I(Y_i \in I_{j,k}) - E[\widetilde{\tau}^{(i)} I(Y_i \in I_{j,k})]\}\right| > C_\lambda \sqrt{m 2^{-j}} \log m\right)
$$
$$
= O(m^{-\lambda}).
$$

Define

$$S_{j,k} = \sum_{i=1}^{m} E\tau^{(i)} I(X_{i-1} \in I_{j,k}) - C_\lambda \sqrt{m 2^{-j}} \log m.$$

Furthermore, we define

$$T_{j,k} = S_{j,k} - \sum_{(u,v):\, I_{u,v} \subset I_{j,k}} S_{u,v}.$$

Then $S_{j,k} = \sum_{(u,v):\, I_{u,v} \subset I_{j,k}} T_{u,v}$. By (39) and (40) we obtain (37) and (38).

(v) *Conclusion for $Z^1_{j,k} - Z^2_{j,k}$.*

By (35)–(38) we obtain that

$$\begin{aligned}
Z^1_{j,k} &= \sum_{(u,v):\, I_{u,v} \subseteq I_{j,k}} W_{u,v}(T_{u,v}) \\
&\quad + \sum_{i:\, 1 \leq i \leq m, X_{i-1} \in I_{j,k}} \sum_{(u,v):\, I_{u,v} \supset I_{j,k}} W_{u,v}(\tau^{(i)}_{u,v}) - W_{u,v}(\tau^{(i-1)}_{u,v})
\end{aligned} \quad (41)$$

and

$$\begin{aligned}
Z^2_{j,k} &= \sum_{(u,v):\, I_{u,v} \subseteq I_{j,k}} W_{u,v}(T_{u,v}) \\
&\quad + \sum_{i:\, 1 \leq i \leq m, Y_i \in I_{j,k}} \sum_{(u,v):\, I_{u,v} \supset I_{j,k}} W_{u,v}(\widetilde{\tau}^{(i)}_{u,v}) - W_{u,v}(\widetilde{\tau}^{(i-1)}_{u,v})
\end{aligned} \quad (42)$$

are satisfied with a probability exceeding $1 - O(m^{-\lambda})$. At this point we see why our particular pairing of the random variables provides a close connection between $Z^1_{j,k}$ and $Z^2_{j,k}$: most of the randomness of $Z^1_{j,k}$ and $Z^2_{j,k}$ is contained in the first terms on the right-hand sides of (41) and (42), respectively. These terms are random, but identical to each other.

To analyze the difference between the right-hand sides of (41) and (42), we compose the pieces $\{W_{u,v}(s), \tau^{(i-1)}_{u,v} \leq s \leq \tau^{(i)}_{u,v}\}$ and $\{W_{u,v}(s), \widetilde{\tau}^{(i-1)}_{u,v} \leq s \leq \widetilde{\tau}^{(i)}_{u,v}\}$ corresponding to intervals $I_{u,v} \supset I_{j,k}$ to Wiener processes. For fixed $i$,

we define

$$W_{j,k}^{\text{res},i}(s) = \begin{cases} W_{j-1,[k/2]}(s + \tau_{j-1,[k/2]}^{(i-1)}) - W_{j-1,[k/2]}(\tau_{j-1,[k/2]}^{(i-1)}), \\ \quad \text{if } 0 \leq s \leq \tau_{j-1,[k/2]}^{(i)} - \tau_{j-1,[k/2]}^{(i-1)}, \\ [W_{j-1,[k/2]}(\tau_{j-1,[k/2]}^{(i)}) - W_{j-1,[k/2]}(\tau_{j-1,[k/2]}^{(i-1)})] + \cdots \\ \quad + [W_{l+1,[k2^{l+1-j}]}(\tau_{l+1,[k2^{l+1-j}]}^{(i)}) \\ \qquad - W_{l+1,[k2^{l+1-j}]}(\tau_{l+1,[k2^{l+1-j}]}^{(i-1)})] \\ \quad + [W_{l,[k2^{l-j}]}(u) - W_{l,[k2^{l-j}]}(\tau_{l,[k2^{l-j}]}^{(i-1)})], \\ \quad \text{if } s = (\tau_{j-1,[k/2]}^{(i)} - \tau_{j-1,[k/2]}^{(i-1)}) + \cdots \\ \qquad + (\tau_{l+1,[k2^{l+1-j}]}^{(i)} - \tau_{l+1,[k2^{l+1-j}]}^{(i-1)}) + (u - \tau_{l,[k2^{l-j}]}^{(i-1)}) \\ \quad \text{with } \tau_{l,[k2^{l-j}]}^{(i-1)} < u \leq \tau_{l,[k2^{l-j}]}^{(i)}. \end{cases}$$

It is clear that $W_{j,k}^{\text{res},i}$ is a Wiener process on the interval $[0, \tau_{j,k}^{\text{res},i}]$, where $\tau_{j,k}^{\text{res},i} = \sum_{(u,v)\,:\,I_{j,k} \subset I_{u,v}} (\tau_{u,v}^{(i)} - \tau_{u,v}^{(i-1)})$.

By the strong Markov property, the remaining parts of the Wiener processes $W_{j,k}$ again form independent Wiener processes, also independent of $\{W_{j,k}^{\text{res},i}, 0 \leq s \leq \tau_{j,k}^{\text{res},i}\}$. Therefore, we can compose all these latter parts to one Wiener process by setting

$$W_{j,k}^{\text{res}}(s) = \begin{cases} W_{j,k}^{\text{res},1}(s), & \text{if } 0 \leq s \leq \tau_{j,k}^{\text{res},1}, \\ W_{j,k}^{\text{res},1}(\tau_{j,k}^{\text{res},1}) + \cdots + W_{j,k}^{\text{res},u-1}(\tau_{j,k}^{\text{res},u-1}) \\ \quad + W_{j,k}^{\text{res},u}(s - \tau_{j,k}^{\text{res},1} - \cdots - \tau_{j,k}^{\text{res},u-1}), \\ \quad \text{if } \tau_{j,k}^{\text{res},1} + \cdots + \tau_{j,k}^{\text{res},u-1} < s \leq \tau_{j,k}^{\text{res},1} + \cdots + \tau_{j,k}^{\text{res},u}. \end{cases}$$

An analogous construction can be made for the $\widetilde{\tau}_{u,v}^{(i)}$, leading to Wiener processes $\widetilde{W}_{j,k}^{\text{res}}$.

If both $\sum_{i=1}^m \tau^{(i)} I(X_{i-1} \in I_{j,k}) \geq S_{j,k}$ and $\sum_{i=1}^m \widetilde{\tau}^{(i)} I(Y_i \in I_{j,k}) \geq S_{j,k}$ are satisfied, then

$$\sum_{i\,:\,1\leq i\leq m, X_{i-1}\in I_{j,k}} \sum_{(u,v)\,:\,I_{j,k}\subset I_{u,v}} W_{u,v}(\tau_{u,v}^{(i)}) - W_{u,v}(\tau_{u,v}^{(i-1)})$$

$$= W_{j,k}^{\text{res}}\left(\sum_{i\,:\,X_{i-1}\in I_{j,k}} \tau^{(i)} - S_{j,k}\right)$$

and

$$\sum_{i\,:\,1\leq i\leq m, Y_i\in I_{j,k}} \sum_{(u,v)\,:\,I_{j,k}\subset I_{u,v}} W_{u,v}(\widetilde{\tau}_{u,v}^{(i)}) - W_{u,v}(\widetilde{\tau}_{u,v}^{(i-1)})$$

$$= \widetilde{W}_{j,k}^{\text{res}}\left(\sum_{i:\,X_{i-1}\in I_{j,k}} \widetilde{\tau}^{(i)} - S_{j,k}\right).$$

Hence, we obtain by (37), (38) and Lemma 1.2.1 in [6], page 29, that, for all $(j,k) \in \mathcal{I}_n$,

$$\operatorname{ess\,sup} P_{f_0}(|Z_{j,k}^1 - Z_{j,k}^2| > r_n''|X_0)$$

$$\leq \operatorname{ess\,sup} P_{f_0}\left(\left|W_{j,k}^{\text{res}}\left(\sum_{i:\,X_{i-1}\in I_{j,k}} \tau^{(i)} - S_{j,k}\right)\right| > r_n''/2 \Big| X_0\right)$$

$$+ \operatorname{ess\,sup} P_{f_0}\left(\left|\widetilde{W}_{j,k}^{\text{res}}\left(\sum_{i:\,X_{i-1}\in I_{j,k}} \widetilde{\tau}^{(i)} - S_{j,k}\right)\right| > r_n''/2 \Big| X_0\right)$$

$$+ O(m^{-\lambda}) = O(m^{-\lambda}),$$

where $r_n'' = C_\lambda (m 2^{-j})^{1/4} \log m$. This completes the proof. □

## 6. Some auxiliary results.

### 6.1. Convergence of stationary distributions.

LEMMA 6.1. *Suppose that $(X_i^f)_{i\geq 0}$ and $(X_i^{f_0})_{i\geq 0}$ are stationary processes obeying (4) with autoregression functions $f$ and $f_0$, respectively, where $|f|, |f_0| \leq M$. Assume that the innovations $(\varepsilon_i)_{i\geq 1}$ are i.i.d. with a density $p$ such that*

$$\rho = \tfrac{1}{2} \sup_{-M \leq x_1 \leq x_2 \leq M} \int_{-\infty}^{\infty} |p(x - x_1) - p(x - x_2)|\, dx < 1.$$

*Then, for the stationary densities $\psi_f$ and $\psi_{f_0}$, it holds that*

$$\int_{-\infty}^{\infty} (\sqrt{\psi_f(x)} - \sqrt{\psi_{f_0}(x)})^2\, dx$$

$$\leq \frac{1}{1-\rho} \sup_{u \in [0, \|f-f_0\|_\infty]} \int_{-\infty}^{\infty} |p(x) - p(x-u)|\, dx.$$

PROOF. We denote by $p^f(x|y) = p(x - f(y))$ and $p^{f_0}(x|y) = p(x - f_0(y))$ the transition densities of the processes $(X_i^f)_{i\geq 0}$ and $(X_i^{f_0})_{i\geq 0}$, respectively. It holds that

$$\int_{-\infty}^{\infty} (\sqrt{\psi_f(x)} - \sqrt{\psi_{f_0}(x)})^2\, dx \leq \int_{-\infty}^{\infty} |\psi_f(x) - \psi_{f_0}(x)|\, dx.$$

Let, for brevity, $\Psi_{f,f_0}(x) = \psi_f(x) - \psi_{f_0}(x)$. From

$$\Psi_{f,f_0}(x) = \int [p^f(x|y) - p^{f_0}(x|y)] \psi_f(y)\, dy + \int p^{f_0}(x|y) \Psi_{f,f_0}(y)\, dy$$

we deduce that

$$\int_{-\infty}^{\infty} |\Psi_{f,f_0}(x)|\,dx$$
$$\leq \int \left[\int |p^f(x|y) - p^{f_0}(x|y)|\psi_f(y)\,dy\right] dx$$
$$\quad + \int \left|\int p^{f_0}(x|y)[\Psi_{f,f_0}(y)]_+ \,dy - \int p^{f_0}(x|y)[\Psi_{f,f_0}(y)]_- \,dy\right| dx$$
$$\leq \int \left[\int |p^f(x|y) - p^{f_0}(x|y)|\,dx\right]\psi_f(y)\,dy$$
$$\quad + \sup_{y_1,y_2}\int |p^{f_0}(x|y_1) - p^{f_0}(x|y_2)|\,dx \int [\Psi_{f,f_0}(y)]_+ \,dy$$
$$\leq \sup_y \int |p^f(x|y) - p^{f_0}(x|y)|\,dx$$
$$\quad + \sup_{y_1,y_2}\int |p^{f_0}(x|y_1) - p^{f_0}(x|y_2)|\,dx \int [\Psi_{f,f_0}(y)]_+ \,dy.$$

The latter implies

$$\int_{-\infty}^{\infty} |\Psi_{f,f_0}(x)|\,dx \leq \sup_{0\leq u \leq \|f-f_0\|_\infty} \int |p(x) - p(x-u)|\,dx$$
$$\quad + \tfrac{1}{2}\sup_{y_1,y_2}\int |p^{f_0}(x|y_1) - p^{f_0}(x|y_2)|\,dx \int |\Psi_{f,f_0}(x)|\,dy.$$

Rearranging the terms we obtain the assertion. $\square$

6.2. *An analogue of Berbee's lemma.*

DEFINITION 6.1. The uniform $\phi$-mixing coefficient between r.v.'s $\xi$ and $\eta$ is defined to be the number

$$\phi(\xi,\eta) = \sup\{|P(A) - P(A|B)| : A \in \sigma(\xi), B \in \sigma(\eta), P(B) \neq 0\}.$$

LEMMA 6.2. *Suppose that $\xi$ and $\eta$ are two random variables with values in $\mathbb{R}^1$ and $\mathbb{R}^d$, respectively, given on the probability space $(\Omega, \mathcal{F}, P)$. Furthermore, we assume that $\xi$ and $\eta$ possess a joint density and that the probability space is rich enough for the definition of a random variable $\Delta \sim Uniform[0,1]$ which is independent of $\xi$ and $\eta$. Then we can construct a random variable $\widetilde{\xi} = \widetilde{\xi}(\xi,\eta,\Delta)$ such that:*

(i) *$\mathcal{L}(\widetilde{\xi}|\eta) = \mathcal{L}(\xi)$ a.s., that is, $\widetilde{\xi}$ is independent of $\eta$ and has the same distribution as $\xi$,*
(ii) *$P(\widetilde{\xi} \neq \xi | \eta) \leq \phi(\xi,\eta)$ a.s.*

PROOF. The idea of the proof is of course closely related to that of the proof of Theorem 2 in [1]. However, since the formulation of our result differs slightly from theirs (they constructed $\widetilde{\xi}$ in such a way that it is close to $\xi$ with a high probability, whereas we are interested in an exact coincidence of $\widetilde{\xi}$ and $\xi$) we decided not to omit this proof.

We denote by $p_\xi(\cdot)$ the marginal density of $\xi$ and by $p_{\xi|\eta}(\cdot|y)$ the conditional density of $\xi$ given $\eta = y$. Define

$$\phi_y = \tfrac{1}{2}\int |p_\xi(x) - p_{\xi|\eta}(x|y)|\,dx = 1 - \int p_\xi(x) \wedge p_{\xi|\eta}(x|y)\,dx.$$

Then $\phi_\eta \leq \phi(\xi, \eta)$ a.s.

If $\phi_\eta = 0$, then $p_\xi(\cdot)$ and $p_{\xi|\eta}(\cdot \mid y)$ coincide and we set $\widetilde{\xi} \equiv \xi$. Otherwise we proceed as follows. With a random variable $\Delta \sim \textit{Uniform}[0,1]$ which is independent of $\xi$ and $\eta$, we set

$$\widetilde{\xi} = \widetilde{\xi}(\xi, \eta, \Delta) = \begin{cases} \xi, & \text{if } p_\xi(\xi) \geq \Delta p_{\xi|\eta}(\xi|\eta), \\ \overline{\xi}, & \text{otherwise}, \end{cases}$$

where $\overline{\xi}$ is an appropriate random variable having the density $[p_\xi(\cdot) - p_\xi(\cdot) \wedge p_{\xi|\eta}(\cdot|\eta)]/\phi_\eta$. The random variable $\overline{\xi}$ is defined via a quantile transform as $G_\eta^{-1}(\frac{\Delta p_{\xi|\eta}(\xi|\eta) - p_\xi(\xi)}{p_{\xi|\eta}(\xi|\eta) - p_\xi(\xi)})$, where

$$G_\eta(y) = \frac{1}{\phi_\eta}\int_{-\infty}^y [p_\xi(x) - p_{\xi|\eta}(x|\eta)]_+\,dx.$$

Now we have

$$P(\widetilde{\xi} = \xi|\eta) = P(p_{\xi|\eta}(\xi|\eta) \wedge p_\xi(\xi) \geq \Delta p_{\xi|\eta}(\xi|\eta)|\eta)$$
$$= E\left(\frac{p_{\xi|\eta}(\xi|\eta) \wedge p_\xi(\xi)}{p_{\xi|\eta}(\xi \mid \eta)} I(p_{\xi|\eta}(\xi|\eta) > 0)\Big|\eta\right)$$
$$= \int p_{\xi|\eta}(x|\eta) \wedge p_\xi(x)\,dx = 1 - \phi_\eta,$$

which implies (ii). Part (i) follows from the construction. □

6.3. *An exponential inequality.* We made use of the following inequality whose proof can be found in [15].

LEMMA 6.3. *Let $\xi$ be a r.v. such that $E\xi = 0$ and $|\xi| \leq a$, for some positive constant a. Then*

$$E\exp(\lambda\xi) \leq \exp(c\lambda^2 E\xi^2), \qquad |\lambda| \leq 1,$$

*where $c = e^a/2$.*

**Acknowledgments.** We thank the referees and an Associate Editor for helpful comments on the paper.

UNIVERSITÉ DE BRETAGNE-SUD
CENTRE DE RECHERCHE YVES COPPENS
CAMPUS DE TOHANNIC
56000 VANNES
FRANCE
E-MAIL: ion.grama@univ-ubs.fr

INSTITUT FÜR STOCHASTIK
FRIEDRICH-SCHILLER-UNIVERSITÄT JENA
ERNST-ABBE-PLATZ 2
07743 JENA
GERMANY
E-MAIL: mneumann@mathematik.uni-jena.de